\documentclass[12pt]{article}
\usepackage{amsfonts}

\usepackage{enumerate}

\usepackage[colorlinks=true, urlcolor=blue,linkcolor=blue, citecolor=blue]{hyperref}

\usepackage{graphics,graphicx,amsmath}

\newcommand{\R}{\mathbb{R}}

\newcommand{\E}{E}
\newcommand{\EE}{\mathrm{E}}

\newcommand{\N}{\mathbb{N}}

\newtheorem{prop}{Proposition}[section]
\newtheorem{lemma}[prop]{Lemma}

\newtheorem{corollary}[prop]{Corollary}
\newtheorem{theorem}[prop]{Theorem}

\def\({\left(}
\def\){\right)}
\newcommand{\cov}{\mathrm{Cov}}
\def\[{\left[}
\def\]{\right]}
\def\real{{\mathord{\mathbb R}}}
\def\inte{{\mathord{\mathbb N}}}
\def\Dom{\mathrm{Dom}}
\def\Var{\mathrm{Var}}
\newcommand{\p}{\mathbb{P}}

\newenvironment{Proof}{\removelastskip\par\medskip
\noindent{\em Proof.} \rm}{\penalty-20\null\hfill$\square$\par\medbreak}

\allowdisplaybreaks

\numberwithin{equation}{section}

\begin{document}
\title{
\huge
Stein approximation for functionals of independent random sequences 
} 

\author{Nicolas Privault\thanks{Division of Mathematical Sciences,
Nanyang Technological University, SPMS-MAS-05-43, 21 Nanyang Link
Singapore 637371. e-mail: {\tt nprivault@ntu.edu.sg}. 
}
\and 
Grzegorz Serafin\thanks{Faculty of Pure and Applied Mathematics, Wroc{\l}aw University of Science and Technology, Ul. Wybrze\.ze Wyspia\'nskiego 27, Wroc{\l}aw, Poland.
e-mail: {\tt grzegorz.serafin@pwr.edu.pl}.}}

\maketitle

\begin{abstract} 
We derive Stein approximation bounds for functionals of uniform random variables, using chaos expansions and the Clark-Ocone representation formula combined with derivation and finite difference operators. This approach covers sums and functionals of both continuous and discrete independent random variables. For random variables admitting a continuous density, it recovers classical distance bounds based on absolute third moments, with better and explicit constants. We also apply this method to multiple stochastic integrals that can be used to represent $U$-statistics, and include linear and quadratic functionals as particular cases. 
\end{abstract}

\baselineskip0.7cm

\section{Introduction}
The Stein and Chen-Stein methods have been developed 
together with the Malliavin calculus 
to derive bounds on the distances between probability laws
on the Wiener and Poisson spaces, cf. \cite{nourdinpeccati},
\cite{utzet2}, \cite{thale} and for discrete Bernoulli sequences,
cf. \cite{nourdin3}, \cite{reichenbachs}, \cite{krokowski}. 
 The results of these works rely on covariance
 representations based on the number
 (or Ornstein-Uhlenbeck) operator $L$ on multiple
 Wiener-Poisson stochastic integrals and its inverse $L^{-1}$.
 Other covariance representations based on
 the Clark-Ocone representation formula have been used in
 \cite{privaulttorrisi3} on the Wiener and Poisson spaces, and in 
 \cite{privaulttorrisi4} for Bernoulli processes.
\\
 
\noindent
This paper focuses on functionals of a countable number
of uniformly distributed random variables, and uses the framework of
\cite{prebub},
cf. also \cite{girunif}, \cite{qtmunif}, 
to derive covariance representations from chaos
expansions in multiple stochastic integrals,
based on a version of the Clark-Ocone formula
with finite difference or derivation operators.
We obtain general bounds on the distance of a random functional
to the Gaussian and gamma distributions using Stein kernels,
see Propositions~\ref{thm:upperdtv}-\ref{thm:dWI}, 
and we also derive specific bounds for multiple stochastic integrals,
see Corollary~\ref{thm:dIZnabla}. 
Other recent approaches to the Stein method for arbitrary
univariate distributions using Stein kernels include \cite{ley}.   
\\
 
\noindent 
When restricted to single stochastic integrals,
our framework applies to sums 
$$
Z_n : = \frac{1}{\sqrt{n}} \sum_{k=1}^n X_k, \qquad n \geq 1, 
$$
 of independent centered random variables 
 $(X_k)_{k\geq 1}$ with variance one.
 This includes the case of discrete random variables 
 and, e.g., sums and polynomials of Bernoulli random variables
 with variable parameters, 
 as a consequence of Proposition~\ref{cor:dwI}, 
 see Proposition~\ref{djkld2}. 
 In addition, this approach yields the general bound
\begin{equation}
\label{dw} 
  d_W( Z_n ,\mathcal N) \leq \frac{2}{n^{3/2}} \sum_{k=1}^n E[|X_k|^3], 
\end{equation} 
where $d_W$ denotes the Wasserstein distance,
see \eqref{ldjklf} below,
which recovers classical results such as the bound 
of Theorem~1.1 in \cite{goldstein}, however with an additional factor two. 
  \\
   
  \noindent
  On the other hand, for random variables
  which admit a continuous density,
  as a consequence of Proposition~\ref{prop:dw+dtv} 
  we find in Proposition~\ref{djkldddas} that 
 \begin{equation}
  \label{klsdf} 
  d_W ( Z_n ,\mathcal N)
  \leq
  \frac{1}{n^{3/2}} 
  \sqrt{ \sum_{k=1}^n \left(
    \int_{-\infty}^\infty \left|
    \frac{1}{F'_k (y) }
 \int_{-\infty}^y x dF_k (x) \right|^2dF_k (y)
 - 1\right)  
  },
\end{equation}
assuming that the cumulative distribution function $F_k$
of $X_k$ admits a non-vanishing density on the support of $X_k$.
This recovers in particular Proposition~3.3 of \cite{privaulttorrisi3}
in the case $n=1$.
For several usual distributions the bound \eqref{klsdf} improves on
  \eqref{dw} which is based on absolute third moments.
  For example in the Gaussian case, \eqref{klsdf} yields
  $d_W ( Z_n ,\mathcal N)=0$
  as expected. For the Gamma and Beta distributions 
  it also yields better constants than \eqref{dw}.
  The bound \eqref{klsdf} may however perform worse
  than \eqref{dw}, or can become infinite if $F'(x)$ becomes
  too close to $0$ on an interval. 
\\
 
   \noindent
    Multiple stochastic integrals with respect to a point process with
   uniform jump times are particularly treated in
   Proposition~\ref{cor:dwI} and \ref{cor:mf}
   and Corollaries~\ref{thm:dIZnabla} and \ref{thm:dIZD},
   with an application to a combinatorial central limit theorem
   for general i.i.d. random sequences in Theorem~\ref{theoremk}.
   \\

   \noindent
   In Section~\ref{s3.1} we consider
   U-statistics, or quadratic functionals of the form 
$$
 Q_n : = \sum_{1 \leq k , l \leq n} a_{k,l} X_k X_l,
$$
 where $(X_k)_{k\geq 1}$ is a sequence of
 normalized independent identically distributed random
 variables, such that $\Var [Q_n]=1$.
 Corollary~\ref{c01} shows that we have the bound
\begin{equation}
\label{jdk} 
  d_W ( Q_n , \mathcal N)\leq
  2 \sqrt{n} L_n^2 \left(
  C +
  \sqrt{
     E[X_1^4] 
      +
    \frac{2}{nL_n^4} \sum_{1 \leq l , p \leq n }
\(\sum_{k=1}^n
      a_{k,l} a_{k,p}
      \)^2
      } 
\right), 
\end{equation} 
where $C = 3 E[X_1^4] + ( E[X_1^4])^2$ and 
$$
L_n^2 : = \max_{1\leq k \leq n} \sum_{l=1}^n a_{k,l}^2, 
$$
which provides a different bound from Theorem~1 in \cite{goetze},
with explicit constants.
In case $a_{2k,2k-1}=1/\sqrt{n}$, the bound \eqref{jdk} yields 
$$
 d_W ( Q_n , \mathcal N)\leq \frac{8 E[X_1^4]}{\sqrt{n}}, 
$$    
   which recovers the known convergence rate in $1/\sqrt{n}$ 
   as on pages 1074-1075 of \cite{goetze}.
   Corollary~\ref{c01.1} provides another bound obtained
   from derivation operators. 
\\

\noindent 
   More generally, our approach 
   applies to functionals of uniformly distributed random
   variables, see Propositions~\ref{prop:dw+dtv} and
   \ref{thm:dWI} which deal respectively with smooth random functionals
   and with multiple stochastic integrals, cf. Proposition~\ref{cor:dwI}.
\\
    
    \noindent
 This paper is organized as follows.
 In Section~\ref{s1} we recall 
 the framework of \cite{prebub} for the construction
 of random functionals of uniform random variables,
 together with the construction of derivation operators
 and the associated stochastic integral (Clark-Ocone)
 decomposition formula. 
 In Section~\ref{s2} we derive Stein approximation
 bounds for the distance of the laws of
 general functionals to the Gaussian and gamma distributions. 
 Section~\ref{s2.1} deals with single stochastic
 integrals which can be used to represent sums of
 independent random variables.
 Section~\ref{s3} treats the general case of multiple stochastic
 integrals, which can be viewed as $U$-statistics.
 Finally, in Section~\ref{s3.1},
 double stochastic integrals are discussed
 with theirs applications to quadratic functionals. 
 In the appendix Section~\ref{s5} we prove a
 multiplication formula for multiple stochastic integrals.
\section{Functionals of uniform random sequences}
\label{s1}
\subsubsection*{Stochastic integrals}
Consider an i.i.d. sequence
$(U_k)_{k\in \inte}$ of uniformly distributed random
variables on the interval $[-1,1]$,
where $\inte := \{0,1,2,\ldots \}$, and 
let the jump process $(Y_t)_{t\in \real_+}$ be defined as 
$$Y_t : =\sum_{k=0}^\infty\mathbf1_{[2k+1+U_k,\infty)}(t),
  \qquad t\in \real_+. 
$$
  We also denote by $({\cal F}_t)_{t\in \real_+}$
  the filtration generated by $(Y_t)_{t\in \real_+}$, and let 
  $$\tilde{\cal F}_t: ={\cal F}_{2k},\qquad 2k\leq t<2k+2,
  \qquad k \in \inte.
  $$
  The compensated stochastic integral
  $$
  \int_0^{\infty}u_t d(Y_t-t/2)
  $$
  with respect to the compensated point process 
  $(Y_t - t/2)_{t\in \real_+}$ can be defined for square-integrable
  $\tilde{\cal F}_t$-adapted processes
  $(u_t)_{t\in \real_+}$ by the isometry relation 
  \begin{equation}
    \label{eq:uv}
    \E\[\int_0^{\infty} \hskip-0.1cm
    u_t d(Y_t-t/2)
    \int_0^{\infty}
    \hskip-0.1cm
    v_t d(Y_t-t/2)\]
=
\E\[\int_0^{\infty} \hskip-0.1cm 
u_t \(v_t - \sum_{k=0}^\infty\mathbf1_{(2k,2k+2]}(t)\int_{2k}^{2k+2}
\hskip-0.1cm
v_r \frac{dr}{2} 
\)\frac{dt}{2}\],
\end{equation} 
  see \cite{prebub},  where $(u_t)_{t\in \real_+}$ and $(v_t)_{t\in \real_+}$ are
square-integrable $\tilde{\cal F}_t$-adapted processes.
This also implies the bound
\begin{equation*}
\E\[\left(
\int_0^{\infty} \hskip-0.1cm
    u_t d(Y_t-t/2)
\right)^2 \]
\leq 
\frac{1}{2}
\E\[\int_0^{\infty} \hskip-0.1cm 
| u_t |^2 dt\],
\end{equation*}
 for $(u_t)_{t\in \real_+}$ a square-integrable $\tilde{\cal F}_t$-adapted process.
\\
  
\noindent
Given $f_1\in L^1(\R_+)\cap L^2(\R_+)$ 
we define the first order stochastic integral 
$$
I_1(f_1)
:= \sum_{k=0}^\infty f_1 (2k+1+U_k )
- \frac{1}{2} \int_0^\infty f_1 (t) dt 
 = \int_0^\infty f_1(t) d(Y_t-t/2).
$$
 Next, given $f_n$ a function which is square integrable
 on $\R_+^n$ and belongs to the space 
$\hat{L}^2(\R_+^n)$ of symmetric functions that
vanish outside of 
$$
\Delta_n : = \bigcup_{
  k_i \not= k_j \geq 0 \atop
  1\leq i \not= j \leq n
}
      [2k_1,2k_1+2]\times \cdots \times [2k_n,2k_n+2],
$$
we define the multiple stochastic integral 
\begin{align} 
\nonumber 
 I_n(f_n)
  &  := \sum_{r=0}^n
    \frac{(-1)^{n-r}}{2^{n-r}} {{n}\choose{r}} 
   \\
  \nonumber
  &  
  \sum_{ k_1\neq\cdots \neq k_r \geq 0}
  \int_0^\infty\cdots \int_0^\infty f_n (2k_1+1+U_{k_1}
  ,\ldots ,2k_r+1+U_{k_r} ,y_1,\ldots ,y_{n-r})dy_1\cdots dy_{n-r}
  \\
  \nonumber
  &=  n!\int_0^\infty \int_0^{t_n}\cdots \int_0^{t_2}f_n(t_1,\ldots ,t_n)d(Y_{t_1}-t_1/2)\cdots d(Y_{t_n}-t_n/2), 
\end{align} 
see \cite{qtmunif} for a construction using a Wick type product,
and \cite{surgailis} for the Poisson point process version.
It is easy to notice, see
\eqref{eq:uv} above and Propositions~4 and 6 of \cite{prebub},
that $(I_n(f_n))_{n\geq 1}$ forms a 
family of mutually orthogonal centered random variables
which satisfy the bound 
\begin{equation}\label{eq:I^2<}
  \E \big[( I_n(f_n))^2\big]\leq n!\left\|f_n\right\|^2_{L^2(\R^n_+,dx/2)},
  \qquad n\geq 1, 
\end{equation} 
which allows us to extend
the definition of $I_n(f_n)$ to all
$f_n \in \hat{L}^2(\R_+^n)$. If in addition we have 
\begin{equation}
\label{ass:int0}
\int_{2k}^{2k+2}f_n (t,*)dt=0, \qquad k \in \inte, 
\end{equation} 
i.e. the function $f_n$ is canonical \cite{surgailisclt},
then the multiple stochastic integral $I_n(f_n)$
can be written as the $U$-statistic 
of order $n$ based on the function $f_n$, i.e. 
\begin{align}
  \label{eq:defI}
 I_n(f_n)&=\sum_{ k_1\neq\cdots \neq k_n \geq 0 }f_n (2k_1+1+U_1,\ldots ,2k_n+1+U_n),
\end{align}
 with the isometry and orthogonality relation
\begin{equation}
\label{eq:EI^2}
\E\[ I_n(f_n) I_m(f_m) \]=
           {\bf 1}_{\{ n = m \} }
           n!\langle f_n , f_m \rangle_{L^2(\R_+,dx/2)^{\circ n}}, 
\end{equation}
 see \cite{prebub} page~589.
Finally, every $X\in L^2(\Omega)$ admits the chaos decomposition 
\begin{equation}\label{eq:decomp}
X=E[X] + \sum_{n=1}^\infty I_n(f_n),
\end{equation}
for some sequence of functions $f_n$ in $\hat{L}^2(\R_+^n)$,
$n\geq 1$, cf. Proposition~7 of \cite{prebub}. 
\subsubsection*{Finite difference operator} 
Consider the finite difference operator $\nabla$
defined on multiple stochastic integrals $X = I_n(f_n)$ as
\begin{equation}\label{eq:nablaaspsi}
  \nabla_tX := X \circ \Psi_t - \frac{1}{2}
  \int_{2\lfloor t/2\rfloor}^{2\lfloor t/2\rfloor+2}X\circ \Psi_s ds,
  \qquad t \in \real_+, 
\end{equation}
where
$$\Psi_t(\omega) : =\(U_1(\omega),\ldots ,U_{\lfloor t/2\rfloor-1}(\omega),t-2\lfloor t/2\rfloor-1,U_{\lfloor t/2\rfloor+1}(\omega),\ldots \),
\qquad t\in \real_+,
$$
cf. Definition~5 and Proposition~10 of \cite{prebub}.
The operator
$\nabla$ does not satisfy the chain rule of derivation, however it
possesses a simple form and it can be easily applied
to multiple stochastic integrals. 
\begin{prop} 
  Given $f_n\in \hat{L}^2(\R_+^n)$, we have 
\begin{align}\label{eq:nablaext}
\nabla_t I_n(f_n)=n I_{n-1}(f_n(t,*))-n\int_{2\lfloor t/2\rfloor}^{2\lfloor t/2\rfloor+2} I_{n-1}(f_n(s,*))ds, \qquad t\in \real_+. 
\end{align}
\end{prop}
\begin{Proof}
  We observe that
\begin{equation}\label{eq:Ipsi}
 I_n(f_n) \circ \Psi_t = I_n(f_n)+n I_{n-1}(f_n(t,*))-n I_{n-1}(f_n(v,*))_{\mid v=2\lfloor t/2\rfloor+1+U_{\lfloor t/2\rfloor}}.
\end{equation}
 Consequently we have 
\begin{eqnarray*}
\lefteqn{
  \frac{1}{2} \int_{2\lfloor t/2\rfloor}^{2\lfloor t/2\rfloor+2}
  I_n(f_n) \circ \Psi_s ds
  }
  \\
  &= & I_n(f_n)+n\int_{2\lfloor t/2\rfloor}^{2\lfloor t/2\rfloor+2} I_{n-1}(f_n(s,*))ds-n I_{n-1}(f_n(v,*))_{\mid v=2\lfloor t/2\rfloor+1+U_{\lfloor t/2\rfloor}},
  \quad
  t\in \real_+,
\end{eqnarray*}
 and applying this to \eqref{eq:nablaaspsi} we obtain the conclusion.
\end{Proof}
 In particular, under the condition \eqref{ass:int0} we have the equality
$$
 \nabla_t I_n(f_n) = n I_{n-1}\(f_n(t,*)\), \qquad
 t\in \real_+, 
$$
 as in Proposition~10 of \cite{prebub}.
 The operator $\nabla$ also admits an adjoint operator $\nabla^*$ given by 
$$
 \nabla^* \left( I_n(g_{n+1}) \right) : = I_{n+1}({\bf 1}_{\Delta_{n+1}}
 \tilde{g}_{n+1}), 
$$ 
where $\tilde{g}_{n+1}$ is the symmetrization of
$g_{n+1} \in \hat{L}^2(\R_+^n)\otimes L^2(\R_+)$
in $n+1$ variables, 
and $\nabla$ is closable with domain
$$
\Dom ( \nabla ) =
\big\{ X\in L^2(\Omega ) : E [ \Vert \nabla X \Vert^2_{L^2(\real_+)} ]
< \infty \big\}
  , 
  $$
and we have the duality relation
\begin{equation}
\label{dr} 
E[\langle \nabla X , u \rangle_{L^2(\R_+,dx/2)} ] 
=
E[ X \nabla^*(u)], \qquad X \in \Dom ( \nabla ),
\end{equation}
 for $u$ in the domain $\Dom ( \nabla^* )$ of $\nabla^*$,
 cf. Proposition~8 of \cite{prebub}. 
The operator $L$ defined on linear combinations of multiple stochastic
integrals as
\begin{equation*}
  L I_n(f_n) : =-\nabla^*\nabla_t I_n(f_n) =-n I_n(f_n),
  \qquad
  f_n\in \hat{L}^2(\R_+^n), 
\end{equation*}
is called the Ornstein-Uhlenbeck operator. By \eqref{eq:decomp}
the operator is well-defined, invertible for centered
$X\in L^2(\Omega)$, and the inverse operator $L^{-1}$ is given by
$$L^{-1} I_n(f_n)  =-\frac1n I_n(f_n),\qquad n\geq1.$$ 
\noindent
 Recall that the operator $\nabla$ satisfies the Clark-Ocone formula
\begin{equation}\label{eq:C-O2}
X=\E\[X\]+\int_0^{\infty}\E\big[\nabla_tX \mid \tilde{{\cal F}}_t\big]d(Y_t-t/2), 
\end{equation}
for $X\in L^2(\Omega )$, see \cite{prebub}, Theorem~2.
This relation is reformulated using the operator $\Psi_t$
in the next proposition. 
\begin{prop} For all $X\in L^2(\Omega )$ we have
\begin{equation}\label{eq:C-O3}
X=\E\[X\]+\int_0^{\infty}\E\big[ X\circ \Psi_t \mid \tilde{{\cal F}}_t\big]d(Y_t-t/2).
\end{equation}
\end{prop}
\begin{Proof}
 Since the integral term in the right hand side of
 \eqref{eq:nablaaspsi} is constant in $t$ on every interval
 of the form $[2k,2k+2)$, $k\in \inte$, we get
\begin{eqnarray*}
\int_0^\infty \E\big[\nabla_tX \mid \tilde{{\cal F}}_t\big]d( Y_t-t/2)
&= & \int_0^\infty \E\[X\circ \Psi_t-\frac{1}{2}
\int_{2\lfloor t/2\rfloor}^{2\lfloor t/2\rfloor+2}X\circ \Psi_s ds \Big| \tilde{{\cal F}}_t\]d( Y_t-t/2)\\
&=&\int_0^\infty \E\big[ X\circ \Psi_t \mid \tilde{{\cal F}}_t\big]d( Y_t-t/2), 
\end{eqnarray*}
 and \eqref{eq:C-O2} ends the proof.
\end{Proof}
In particular, it follows from the Clark-Ocone formula \eqref{eq:C-O2}
that 
\begin{align}\label{eq:COI}
  \int_0^\infty\E\big[ I_{n-1}(f_n(t,*))
  \mid \tilde{{\cal F}}_t\big] d(Y_t-t/2)=\frac1n I_n(f_n), 
\end{align}
since the integral term in the right hand side of \eqref{eq:nablaext}
is constant in $t$ on every interval of the
form $[2k,2k+2)$, $k\in \inte$.
\subsubsection*{Derivation operator}
\noindent
 Given $X$ a random variable of the form
$$
X=f(U_0,\ldots ,U_n),
\qquad
f\in {\cal C}_b^1([-1,1]^{n+1}),
$$
 we consider the gradient $D_t$ defined as 
$$D_tX :=\sum_{k=1}^{n}\partial_kf(U_0,\ldots ,U_n)
\left(
(1-U_k)\mathbf1_{(2k,2k+1+U_k]}(t)-(1+U_k)\mathbf1_{(2k+1+U_k,2k+2]}(t)
 \right),
 $$ 
 cf. Definition~3 of \cite{prebub}.
 By Proposition~5 of \cite{prebub} the gradient $D$ is closable, 
 and its closed domain is denoted by $\Dom (D)$.
    For any $X\in\Dom ( D )$ and $\phi\in {\cal C}^1_b(\real )$
    we have $\phi(X)\in\Dom ( D )$,
    and the operator $D$ satisfies the chain rule of derivation
\begin{equation}
\label{eq:derule2}
 D_t\phi(X)=\phi'(X) D_t X,
 \qquad
 X\in\Dom ( D ),
\end{equation}
 for all $\phi\in {\cal C}^1_b(\real )$. 
 The gradient operator
$$
D :\Dom ( D )\subset L^2(\Omega )\longrightarrow
L^2(\Omega\times \R_+ )
$$
with domain $\Dom ( D )$,
defined by $D X=( D_t X)_{t\in \real_+}$
satisfies the following Clark-Ocone representation formula,
see Theorem~2 of \cite{prebub}. 
\begin{prop}For $X\in L^2(\Omega )$ we have
  \begin{equation}\label{eq:C-O}
  X=\E\[X\]+\int_0^{\infty}\E\big[ D_tX \mid \tilde{{\cal F}}_t\big]d(Y_t-t/2). 
\end{equation}
\end{prop}
\subsubsection*{Covariance identities} 
 From \eqref{eq:derule2} the gradient operator $D$ satisfies the 
 following covariance identity, see e.g.
 Proposition~3.4.1 in \cite{privaultbk2}, p.~121.
\begin{lemma}\label{lem:cov}
 Let $X,Y \in \Dom ( D )$. We have
 $$\cov(X,Y)=\frac12\E\[\int_0^\infty\E\big[ D_tX
 \mid \tilde{{\cal F}}_t\big] D_tY\,dt\].
$$
\end{lemma}
\begin{Proof}
By \eqref{eq:uv} and \eqref{eq:C-O} we have 
\begin{align*}
\cov(X,Y)&=\E\[(X-\E\[X\])(Y-\E\[Y\])\]\\
&=\E\[\int_0^{\infty}\E\big[ D_tX\mid \tilde{{\cal F}}_t\big]d(Y_t-t/2)\int_0^{\infty}\E\big[ D_tY\mid \tilde{{\cal F}}_t\big] d(Y_t-t/2)\]\\\nonumber
&=\frac{1}{2}
\E\[\int_0^{\infty}\E\big[ D_tY\mid \tilde{{\cal F}}_t\big]\(\E\big[ D_tX\mid \tilde{{\cal F}}_t\big]-\Phi_t(X)\)dt
\],
\end{align*}
where
$$
\Phi_t(X) : =
\frac{1}{2}
\sum_{k=0}^\infty\mathbf1_{(2k,2k+2]}(t)\int_{2k}^{2k+2}\E\big[ D_rX \mid \tilde{{\cal F}}_r\big]dr,
  \qquad
  t\in\real_+.
$$
  By the independence between
  $\tilde{\cal F}_{2k}$ and $(U_k,\ldots ,U_{n})$ we get
\begin{align*}
  &     \Phi_t(X)=\frac12\sum_{k=0}^n\mathbf1_{(2k,2k+2]}(t)
      \\
    &  
    \times \int_{2k}^{2k+2}\E \big[ \partial_kf(U_1,\ldots ,U_n)\((1-U_k)\mathbf1_{(2k,2k+1+U_k]}(r) \right.
    \left.\left.-(1+U_k)\mathbf1_{(2k+1+U_k,2k+2]}(r)\) \big| \tilde{\cal F}_{2k}\]dr\\
    &=  \frac12\sum_{k=0}^n\mathbf1_{(2k,2k+2]}(t)\int_{2k}^{2k+2}\E\[\partial_kf(y_1,\ldots ,y_{k-1},U_k,\ldots ,U_n)
\right.
        \\
        & 
\qquad        \left.
        \times \((1-U_k)\mathbf1_{(2k,2k+1+U_k]}(r)
          -(1+U_k)\mathbf1_{(2k+1+U_k,2k+2]}(r)\)\]_{\mid(y_1,\ldots ,y_{k-1})=(U_1,\ldots ,U_{k-1})}dr\\
&=  \frac12\sum_{k=0}^n\mathbf1_{(2k,2k+2]}(t)\E\[\partial_kf(y_1,\ldots ,y_{k-1},U_k,\ldots ,U_n)\right.\\
            &  
              \qquad \qquad\qquad \qquad
              \times
            \left. \big((1-U_k)(1+U_k)-(1-U_k)(1+U_k)\big)\]_{\mid (y_1,\ldots ,y_{k-1})=(U_1,\ldots ,U_{k-1})}\\
&=  0.
\end{align*}
We conclude that 
\begin{align*}
\cov(X,Y)&=\frac12\E\[\int_0^{\infty}\E\big[ D_tY\mid \tilde{{\cal F}}_t\big]\E\big[ D_tX \mid \tilde{{\cal F}}_t\big]dt\]\\
&=\frac12\E\[\int_0^{\infty}\E\big[
  \E\big[ D_tX\mid \tilde{{\cal F}}_t\big] D_tY\mid \tilde{{\cal F}}_t\big]dt\]\\
&=\frac12\E\[\int_0^{\infty}\E\big[ D_tX\mid \tilde{{\cal F}}_t\big] D_tY\,dt\].
\end{align*}
\end{Proof}
\noindent
As a consequence of Lemma~\ref{lem:cov} we have
the inequality 
\begin{equation} 
  \label{inequ}
  \frac{1}{2}
  \EE[\langle
 {D}_\cdot X
 ,
 \EE[ D_\cdot X \mid \tilde{\cal F}_\cdot ]\rangle_{L^2(\R_+ )}
  ] = \Var [ X ] \leq
  \|X\|_{L^2(\Omega )}^2.
\end{equation} 
Using the operator $\nabla$ and the Clark-Ocone formula
\eqref{eq:C-O2}-\eqref{eq:C-O3}
we can also obtain the covariance identity 
  $$\cov(X,Y)=\frac12\E\[\int_0^\infty\E\big[ \nabla_tX
  \mid \tilde{{\cal F}}_t\big] \nabla_tY\,dt\]
$$
from \eqref{eq:uv} and \eqref{eq:nablaaspsi}
as in the proof of Lemma~\ref{lem:cov}.
\subsubsection*{Stein kernel} 
The next proposition shows that the Stein kernel 
$\varphi_X$ defined in \eqref{varphi}
is a Stein kernel in the sense of Definition~(2.1)
in \cite{lnp}.
\begin{prop}
 Let $X\in\Dom (D)$ be such that $\mathrm{E}[X]=0$.
 The Stein kernel 
\begin{equation} 
\label{varphi}
 \varphi_X(z):=
 \frac{1}{2}
 \E \left[\int_0^{\infty} D_t X \E \big[ D_t
X\mid \tilde{\cal F}_t\big] dt\,\Big|\,X=z\right],
 \qquad z\in\R, 
\end{equation} 
satisfies
\begin{equation}
  \label{id} 
  \mathrm{Cov}(X,\phi(X)) = \E [\phi'(X)\varphi_X(X)], 
\end{equation} 
 for any $\phi\in\mathcal{C}_b^1(\mathbb{R})$. 
\end{prop} 
\begin{Proof}
  We note that by Lemma~\ref{lem:cov}
  and Jensen's inequality we have 
$$
 \|\varphi_X(X)\|_{L^1(\Omega )}\leq
 \|\varphi_X(X)\|_{L^2(\Omega )} \leq \E[X^2] = \sqrt{
   \E \left[\int_0^{\infty}| D_t X|^2 \frac{dt}{2} \right]}
 <\infty,
$$
 and, for any $\phi\in\mathcal{C}_b^1(\mathbb{R})$, 
\begin{align}
  \nonumber
  \mathrm{Cov}(X,\phi(X))&= \frac{1}{2} \E \left[\int_{0}^{\infty}
  \E [D_t X\,|\,\mathcal{F}_t]\,D_t \phi(X)dt\right]
\\
\nonumber
&= \frac{1}{2}
\E \left[\phi'(X)\int_{0}^{\infty}D_t X\,\E [D_t X\,|\,\mathcal{F}_t]dt\right]
\\
\nonumber
&=\frac{1}{2}
\E \left[\E \left[\phi'(X)\int_{0}^{\infty}D_t X\,\E [D_t X\,|\,\mathcal{F}_t]dt\,\Big|\,X\right]\right]\\
\label{ibp1}
&=\E [\phi'(X)\varphi_X(X)]. 
\end{align}
\end{Proof} 
In particular, \eqref{ibp1} shows that we have 
$$
\E [ \varphi_X(X)] = \Var [ X ],
\qquad
X\in \Dom (D).
$$ 
In the sequel we will also use the identity
\begin{equation}
   \label{djkf} 
 \varphi_{X_k} ( y )= - \frac{1}{F_k'( y)} \int_{-\infty}^y x dF_k (x),
\end{equation} 
see Relation~(3.17) in \cite{viens}. Next, we review some
examples of Stein kernels. 
\begin{enumerate}[]
\item {\bf Gaussian case.} 
  The Stein kernel
  of ${X_1} \simeq {\cal N}(0,\sigma^2)$ with the
 Gaussian cumulative distribution function $F(x)$
  is given by 
$$ 
 \varphi_{X_1} ( y ) = 
 - \frac{1}{F' (y)} \int_{-\infty}^y x dF (x)
 = 
 \sigma^2, \qquad y \in \real. 
$$ 
\item {\bf Gamma case.}
\noindent 
 When ${X_1}$ has the centered gamma distribution with shape parameter $s>0$
 and density function 
$$
F_s '(x) =
  \frac{(x+s)^{s-1}}{\Gamma (s)} e^{-(x+s)}
  \qquad
  x\in [-s,\infty ), \qquad k \geq 1, 
$$
  we have $E[|{X_1}-s|]=2s^{s} e^{-s}$, hence 
  the Stein kernel of $X_1$ is 
  $$ 
 \varphi_{X_1} ( y ) = 
 - \frac{1}{F'_s (y)} \int_{-s}^y x dF_s (x)
 = y+s, \qquad y \in \real. 
$$ 
\item {\bf Beta case.}
\noindent
When ${X_1}$ has the centered
Beta$(\alpha, 1)$ distribution, $\alpha >0$, we have
$$
 F_\alpha (x) = \left( \frac{\alpha}{\alpha + 1 } + x 
 \right)^\alpha,
 \qquad
 x\in \left[ 
 - 
 \frac{\alpha}{\alpha + 1}
 ,
 \frac{1}{\alpha + 1}
 \right], 
$$
 and the Stein kernel of $X_1$ is 
 \begin{equation}
   \label{djkd} 
 \varphi_{X_1} ( y ) = 
 - \frac{1}{F'_\alpha  (y)} \int_{-\alpha / (\alpha + 1 ) }^y x dF_\alpha  (x)
 = \frac{1}{\alpha + 1}
 \left( \frac{\alpha}{\alpha + 1} + y \right)
 \left( \frac{1}{\alpha + 1} - y \right), \quad y \in \real. 
 \end{equation}
\item {\bf Single stochastic integrals.}
Such integrals can be used to represent
the sum $Z_n$ of independent centered random variables
$(X_k)_{k\geq 1}$ as 
\begin{equation}
  \label{zn} 
  Z_n = \sum_{k=1}^n X_k
   = I_1 \big( f_1 {\bf 1}_{[0,2n]} \big)
  , 
\end{equation} 
where 
\begin{equation}
\label{djkldd} 
f_1(t):= \sum_{k=0}^\infty F_k^{-1} \left( \frac{t}{2} -k \right) \mathbf1_{[2k,2k+2)}(t),
\end{equation} 
satisfies
$\int_{2k}^{2k+2}f_1 (t)dt=0$, $k \in \inte$,
and
$$F_k^{-1} (t) : = \inf \{ s \in \real_+ \ : \ F_X(s) \geq t \},
\qquad
t\in [0,1],
$$
 is the right-continuous inverse of
 the cumulative distribution function $F_k$ of $X_k$, $k\geq 1$.
\end{enumerate} 
 In the sequel we let $C^1_{\Box}(\R_+)$ denote the set of
 functions which are ${\cal C}^1$ on every interval of the
 form $(2k,2k+2)$, $k\in \N$.
 The next lemma can be useful when computing the Stein
 kernel of single stochastic integrals according to \eqref{varphi},
 see Propositions~\ref{djkdddfds} and \ref{djkldddas} below. 
\begin{lemma}
\label{alm} 
Assume that
$\displaystyle Z_n = I_1 \big( f_1 {\bf 1}_{[0,2n]} \big)
= \sum_{k=1}^n X_k
$ 
 belongs to $\Dom (D)$, $n\geq 1$.
We have 
          $$
     \langle D_\cdot Z_n ,
    \E\big[ D_\cdot Z_n \mid \tilde{{\cal F}_\cdot }\big] \rangle_{L^2(\R_+ )}
    = - 2 I_1(
      \varphi_{X_{1+[\cdot /2]}} ( f_1( \cdot ) ) ) + E[Z_n^2]. 
$$ 
        \end{lemma}
      \begin{Proof}
 We note that for $f_1 \in C^1_{\Box}(\R_+) \cap L^2(\R_+)$, we have
 $$
 D_t I_1(f_1)
   = 
   \sum_{k=0}^\infty
   \left(
  (1-U_k )\mathbf1_{(2k,2k+1+U_k]}(t)-(1+U_k)\mathbf1_{(2k+1+U_k,2k+2]}(t)
      \right) f_1' (2k+1+U_k). 
$$ 
 Next, by Proposition~10 and Lemma~1 in \cite{prebub} we get
\begin{align*}
  \E\big[ D_t I_1(f_1) \mid \tilde{{\cal F}}_{t}\big]=\E\big[\nabla_t I_1(f_1) \mid \tilde{{\cal F}}_{t}\big]= f_1(t),
  \qquad t\in \real_+, 
\end{align*}
 hence by \eqref{ass:int0} we have 
\begin{eqnarray} 
  \nonumber
  \lefteqn{
    \langle D_\cdot I_1(f_1) ,
    \E\big[ D_\cdot I_1(f_1) \mid \tilde{{\cal F}_\cdot }\big] \rangle_{L^2(\R_+ )}
}
\\
  \nonumber
  &= & \int_0^\infty \sum_{k=0}^\infty
  \left(
  (1-U_k)\mathbf1_{(2k,2k+1+U_k]}(s)-(1+U_k)\mathbf1_{(2k+1+U_k,2k+2]}(s)
      \right) 
f_1' (2k+1+U_k ) f_1(s)ds
      \\
      \nonumber
      &=&
      \int_0^\infty \sum_{k=0}^\infty 
      \left( \mathbf1_{(2k,2k+1+U_k]}(s)-\mathbf1_{(2k+1+U_k,2k+2]}(s)\right)
 f_1' (2k+1+ U_k ) f_1(s) ds
          \\
          \nonumber
          &=& 2\int_0^\infty
          \sum_{k=0}^\infty
          \left(
          \mathbf1_{(2k,2k+1+U_k]}(s) f_1' (2k+1+U_k)
          f_1(s)\right)
          ds
          \\
          \nonumber
          &=& 2\sum_{k=0}^\infty f_1' (2k+1+U_{k}) \int_0^{2k+1+U_{k}}f_1(s)
          ds
          \\
          \nonumber
          &=& 2\int_0^\infty
          f_1' (t) \int_0^t f_1(s) ds d(Y_t - t/2)
          +\int_0^\infty f_1' (t) \int_0^t f_1(s) dsdt
             \\
\nonumber 
             &= & 2\int_0^\infty
             f_1' (s) \int_0^t f_1(s) dsd(Y_t-t/2)+
             \int_0^\infty |f_1 (t)|^2 dt. 
\end{eqnarray}
On the other hand, by 
\eqref{eq:uv} and
 \eqref{djkf}, see (3.17) in \cite{viens}, we have 
\begin{eqnarray*}
     f'_1(x) \int_0^{x}f_1(t)dt 
  & = &  
 \frac{1}{2} \sum_{k=0}^\infty \frac{1}{F_k'(F_k^{-1} ( (x-2k)/2 ) ) } 
    \int_{2k}^x
    F_k^{-1} \left( \frac{t}{2} - k \right) dt
   \mathbf1_{[2k,2k+2)}(x) 
      \\
  \nonumber 
  & = &
  \sum_{k=0}^\infty
  \mathbf1_{[2k,2k+2)}(x)
    \frac{1}{F_k'(F_k^{-1} ( (x-2k)/2 ) ) }
      \int_{-\infty}^{F_k^{-1} ((x-2k)/2)}
    t dF_k (t) 
      \\
  \nonumber 
  & = &
  -
  \sum_{k=0}^\infty
 \varphi_{X_k} ( F_k^{-1} ( (x-2k)/2 ) 
 \mathbf1_{[2k,2k+2)}(x) 
      \\
  \nonumber 
  & = &
  -
  \sum_{k=0}^\infty
 \varphi_{X_k} ( f_1(x) ) \mathbf1_{[2k,2k+2)}(x) 
      \\
  \nonumber 
  & = &
  -
 \varphi_{X_{1+[x/2]}} ( f_1(x) ) \mathbf1_{[2k,2k+2)}(x) 
, 
\end{eqnarray*}
where we used the identity \eqref{djkf}.
\end{Proof}
\subsubsection*{Density representation and bounds}
Working along the lines of the proof of
Theorem 3.1 in \cite{viens}
by replacing (3.15) therein with
\eqref{ibp1} above we can derive the
following result,
where $\mathrm{Supp}(f)$ denotes the support of the function $f$.
\begin{prop}
\label{prop:denspoiss}
Let $X\in\Dom (D)$ be such that $\mathrm{E}[X]=0$.
 The law of $X$ has a density $p_X$
 with respect to the Lebesgue measure if and only if the
 Stein kernel $\varphi_X$ defined in \eqref{varphi}
 satisfies $\varphi_X(X)>0$ a.s.
 In this case $\mathrm{Supp}(p_X)$ is
 a closed interval of $\real$ containing $0$
 and we have
 \begin{equation*}
    p_X(z)=\frac{\mathrm{E}[|X|]}{2\varphi_X(z)}\exp\left(
 -\int_{0}^z
 \frac{u}{\varphi_X(u)}\,\mathrm{d}u\right),\quad\text{a.e. $z\in\mathrm{Supp}(p_X)$.}
\end{equation*}
\end{prop} 
 As a consequence of Proposition~\ref{prop:denspoiss}
 we get the following result on density bounds
 as in Corollary~3.5 of \cite{viens}. 
 \begin{prop}
   Let $X \in\Dom (D)$ be a centered random variable such that
 \begin{equation*}
   0< c \leq\int_0^\infty D_s X\E [D_s X\,|\,\mathcal{F}_s]ds\leq C \quad a.s.,
\end{equation*}
 where $C , c >0$ are positive constants.
 Then 
 the density $p_X$ satisfies
 \begin{equation*}
   \frac{\E [|X|]}{2C}\,\exp\left(-\frac{z^2}{2c}\right)\leq p_X(z)\leq\frac{\E [|X|]}{2c}\,
\exp\left(-\frac{z^2}{2C}\right),\quad a.e. \ z\in\R, 
\end{equation*}
and the tail probabilities satisfy
\begin{equation*}
  P(X\geq x)\leq\exp\left(-\frac{x^2}{2C}\right)\quad\text{and}\quad P(X\leq
-x)\leq\exp\left(-\frac{x^2}{2C}\right),\qquad x>0.
\end{equation*}
 \end{prop}

 \section{Stein approximation bounds}
\label{s2}
 The total variation distance
 between two real-valued random variables
 $X$ and $Y$ is defined by
 $$d_{TV}(X,Y)=\sup_{A\in\mathcal B(\R)}\left|\p\(X\in A\)-\p\(Y\in A\)\right|,
 $$
 where $\mathcal{B}(\R)$ denotes the Borel subsets of $\R$. 
 The Wasserstein distance between
 the laws of $X$ and $Y$ is defined by
$$
 d_W (X,Y):
 =\sup_{h\in\mathrm{Lip}(1)}|\E [h(X)]-\E [h(Y)]|,
$$
 where $\mathrm{Lip}(1)$ is the class of real-valued
 Lipschitz functions with
 Lipschitz constant less than or equal to $1$. 
 \\

\noindent
In the following propositions
we derive bounds for the 
 Wasserstein and total variation distances between the
 normal distribution and the distribution of 
 a given random variable $X\in\Dom ( D)$. 
 Recall that by Stein's lemma,
 cf. \cite{stein}, \cite{nourdin},
 for any continuous function
 $h : \R \longrightarrow [0,1]$ 
 the Stein equation
 \begin{equation*}
   h(x)-\E\[h(X)\]={f}_h'(x)-xf_h(x),
   \end{equation*} 
 where $X\sim \mathcal N$,
 admits a solution $f_h(x)$ that
 satisfies the bound $|f_h^{'}(x)|\leq 2$. 
 In the sequel we denote by
$$
 \mathcal{T} :=\big\{h\in\mathcal{C}_b^2(\R) \ : \
 \|h'\|_\infty \leq 1, \ \|h''\|_\infty \leq 2
 \big\}
$$ 
 the space
 of twice differentiable functions whose first derivative is bounded by $1$
 and whose second derivative is bounded by $2$. 
 For the gamma approximation we will use the distance
\begin{equation}
\nonumber 
d_{\mathcal{H}}(X,Y):=\sup_{h\in\mathcal{H}}|\E [h(X)]-\E [h(Y)]|,
\end{equation}
 where
$$
\mathcal{H}:=\big\{h\in\mathcal{C}_b^2(\R) \ : \ \max\{\|h\|_\infty,\|h'\|_\infty,\|h''\|_\infty\}\leq 1\big\}.
$$ 
\subsubsection*{Derivation operator bounds}
In the next Proposition~\ref{thm:upperdtv} we derive
a Stein bound using 
the Stein kernel $\varphi_X(z)$ defined in 
\eqref{varphi},
see also Proposition~3.3 of \cite{privaulttorrisi3}
for a bound using a different probabilistic representation of the
Stein kernel. 
Here we denote by 
 $\Gamma(\nu/2)$ a random variable distributed according to the gamma law
 with parameters $(\nu/2,1)$, $\nu>0$.
 We also let $\langle\cdot,\cdot\rangle$ denote the usual inner product
 $\langle\cdot,\cdot\rangle_{L^2(\R_+)}$ on $L^2(\R_+)$.
\begin{prop}\label{thm:upperdtv}
  For any $X\in\Dom ( D )$ such that $\E [X]=0$, we have 
$$ 
 d_W (X,\mathcal{N})
 \leq 
 \mathrm{E}[|1-\varphi_X (X)|]
 \leq 
 | 1 - E[X^2] |
 +
 \Vert
 \varphi_X(X)
 -
 E[ \varphi_X(X)] 
 ]
 \Vert_{L^2(\Omega )}, 
$$ 
where the Stein kernel $\varphi_X$ is defined in \eqref{varphi},
and
$$ 
 d_{TV} (X,\mathcal{N})
 \leq 
 2\mathrm{E}[|1-\varphi_X(X)|]
 \leq 
 2 | 1 - E[X^2] |
 +
 2 \Vert
 \varphi_X(X)
 -
 E[ \varphi_X(X)]
 ]
 \Vert_{L^2(\Omega )}. 
$$ 
 If moreover $X$ is a.s. $(-\nu,\infty)$-valued then we have
$$
d_{\mathcal{H}} (X\, , \, \Gamma_\nu)
\leq\mathrm{E}[|2(X+\nu)-\varphi_X(X)|]
\leq
\Vert 2(X+\nu)-E[X^2]\Vert_{L^2(\Omega )}
+ 
\Vert \varphi_X(X) - E[ \varphi_X(X)] \Vert_{L^2(\Omega )}.
$$
\end{prop}
\begin{Proof}
  We focus on the first inequalities,
  as the second inequalities follow
  from the triangle inequality and Jensen's inequality, 
  and the identity $E[ \varphi_X(X)] = E[X^2]$
  that follows from Lemma~\ref{lem:cov}. 
  \\
  \noindent
 $(i)$ By Lemma \ref{lem:cov} we have 
\begin{eqnarray} 
  \nonumber
  \E\[Xf_h(X)\] & = & \frac12\E\[\int_0^\infty\E\big[ D_tX \mid \tilde{{\cal F}}_t\big]
  D_tf_h(X)\,dt\]
  \\
  \label{cov}
  & = & \frac12\E\[f_h'(X)\int_0^\infty\E\big[ D_tX \mid \tilde{{\cal F}}_t\big] D_tX\,dt\]. 
\end{eqnarray} 
 Hence, using the bound $(2.33)$ in \cite{utzet2} and \eqref{cov}, we get 
\begin{align}
  \label{fjhdfs} 
  d_W(X , \mathcal{N})&\leq\sup_{\phi\in\mathcal{T}}|\EE[\phi'(X)-X\phi (X) ]|
  \\
  \nonumber
  &=\sup_{\phi\in\mathcal{T}} \left|
\EE\left[
  \phi'(X)\left(1-
  \frac{1}{2}
  \langle
 D_{\cdot} X
 ,
 \EE\big[ D_{\cdot} X \mid \tilde{\cal F}_\cdot \big]
 \rangle
 \right) \right]
\right|
  \\
  \nonumber
  &=\sup_{\phi\in\mathcal{T}} \left|
\EE\left[
  \phi'(X)\left(1-
  \varphi_X(X) 
 \right) \right]
\right|
\\
\nonumber
& \leq\EE\left[
  \left| 1- \varphi_X(X) \right| \right].
\end{align}
\noindent 
 $(ii)$
 By the covariance identity \eqref{cov} we have
\begin{align}
|\E [h(X)]-\E [h(\mathcal
  N)]|
&=\left|\E \left[
  f_{h}'(X)
 \left( 1 - 
   \frac{1}{2}
  \langle
 D_{\cdot} X
 ,
 \EE\big[ D_{\cdot} X \mid \tilde{\cal F}_\cdot \big]
 \rangle
 \right)\right]\right|
\nonumber\\
&=\left|\E \left[
  f_{h}'(X)
 \left( 1 - \varphi_X(X) \right)\right]\right|
\nonumber\\
&\leq 2\E \left[\left|1- \varphi_X(X) \right|\right],
\nonumber
\end{align}
and this bound can be extended to $h={\bf 1}_C$
for any $C\in\mathcal{B}_b(\R)$ by the same approximation argument
as in the proof of e.g. Theorem~2.1 of \cite{privaulttorrisi3}. 
\\
\noindent
$(iii)$
 Given $h\in\mathcal{H}$ a twice differentiable
 function bounded above by $1$ we choose $c>0$ and $a<1/2$ such that 
$$
 |h(x)|\leq c\mathrm{e}^{ax},
 \qquad
 x>-\nu. 
$$
 By e.g. Lemma~1.3-$(ii)$ of \cite{nourdinpeccati},
 letting $\Gamma_\nu :=2\Gamma(\nu/2)-\nu$,
 the functional equation
\begin{equation*}
 2(x+\nu)f'(x)=xf(x)+h(x)-\E [h(\Gamma_\nu)],\quad\text{$x>-\nu$},
\end{equation*}
 has a solution $f_h$ which is bounded and
 differentiable on $(-\nu,\infty)$,
 and such that
$$
 \|f_h\|_\infty\leq2\|h'\|_\infty
 \quad
 \mbox{and}
 \quad
 \|f_h'\|_\infty\leq\|h''\|_\infty.
$$
 By the covariance identity \eqref{cov}
on $\mathcal{C}_b^1(\R)$ for the centered random variable $X$ 
 we have
\begin{align}
|\mathrm{E}[h(X)]-\mathrm{E}[h(\Gamma_\nu)]|&=|\mathrm{E}[(2(X+\nu)f_h'(X)-X f_{h}(X))]|\nonumber\\
&= \left| \mathrm{E}\left[f_{h}'(X) \left(2(X+\nu)
 - \langle
 D_{\cdot} X
 ,
 \EE \big[ D_{\cdot} X \mid \tilde{\cal F}_\cdot \big]
 \rangle
 \right)\right] \right|
\nonumber\\
&=\left|
\mathrm{E}\left[f_{h}'(X) \left(2(X+\nu)
 - \varphi_X(X) \right)\right] \right|
\nonumber\\
&\leq \|h''\|_\infty\mathrm{E}[|2(X+\nu)
   - \varphi_X(X) |].\nonumber
\end{align}
 The claim follows by taking the supremum over all
 functions $h\in\mathcal{H}$.
\end{Proof} 
 As a consequence of Proposition~\ref{thm:upperdtv},
for any $X\in\Dom ( D )$ such that $\E [X]=0$, we have 
\begin{eqnarray}
  \nonumber 
 d_W (X,\mathcal{N})
& \leq &
 | 1 - E[X^2] |
 +
 \sqrt{E[ ( \varphi_X(X) - E[X^2] )^2]}
 \\
  \nonumber 
 & = &
| 1 - E[X^2] |
 +
 \sqrt{E[ ( \varphi_X(X))^2 - 2 \varphi_X(X) E[X^2] + ( E[X^2] )^2]}
 \\
 \label{sqrt}
 & = &
| 1 - E[X^2] |
 +
 \sqrt{E[ ( \varphi_X(X))^2 ] - ( E[X^2] )^2}
, 
\end{eqnarray} 
 and 
$$ 
 d_{TV} (X,\mathcal{N})
 \leq 
2 | 1 - E[X^2] |
 +
 2\sqrt{E[ ( \varphi_X(X))^2 ] - ( E[X^2] )^2 }
.
$$ 
 Similarly, Proposition~\ref{thm:upperdtv} implies
 the following corollary which applies in particular to
 smooth functionals $X\in \Dom (D)$. 
 \begin{prop}
   \label{prop:dw+dtv}
 For any $X\in\Dom ( D )$ such that $\E [X]=0$, we have
\begin{eqnarray*}
 d_W (X,\mathcal{N})
 & \leq & 
 \frac{1}{2}
 \Vert
 2 - \langle
 D_{\cdot} X
 ,
 \EE\big[ D_{\cdot} X \mid \tilde{\cal F}_\cdot \big]
 \rangle
 \Vert_{L^2(\Omega )}
\\
 & \leq &
 | 1
 -
 E [ X^2 ] 
 |
 +
 \frac{1}{2}
 \Vert
 \langle
 D_\cdot X
 ,
 \EE\big[ D_\cdot X \mid \tilde{\cal F}_\cdot \big]
 \rangle
 -
 \EE[ \langle
 D_\cdot X
 ,
 \EE\big[ D_\cdot X \mid \tilde{\cal F}_\cdot \big]
 \rangle
 ]
 \Vert_{L^2(\Omega )}, 
\end{eqnarray*}
and
\begin{eqnarray*}
 d_{TV} (X,\mathcal{N})
 & \leq &
 \Vert
 2 
 -
 \langle
 D_{\cdot} X
 ,
 \EE\big[ D_{\cdot} X \mid \tilde{\cal F}_\cdot \big]
 \rangle
 \Vert_{L^2(\Omega )}
\\
 & \leq &
 2
 |
 1
 -
 E[ X^2 ] 
 |
 +
 \Vert
 \langle
 D_\cdot X
 ,
 \EE\big[ D_\cdot X \mid \tilde{\cal F}_\cdot \big]
 \rangle
 -
 \EE[ \langle
 D_\cdot X
 ,
 \EE\big[ D_\cdot X \mid \tilde{\cal F}_\cdot \big]
 \rangle
 ]
 \Vert_{L^2(\Omega )}. 
\end{eqnarray*}
 For any a.s. $(-\nu,\infty)$-valued 
 $X\in\Dom (D)$ such that $\E [X]=0$, we have 
\begin{eqnarray*}
  \lefteqn{
    d_{\mathcal{H}} (X,\Gamma_\nu)
 \leq 
 \Vert
 2(X+\nu)
 -
 \langle
 D_{\cdot} X
 ,
 \EE \big[ D_{\cdot} X \mid \tilde{\cal F}_\cdot \big]
 \rangle
 \Vert_{L^2(\Omega )}
  }
  \\
 & \leq &
 \Vert
 2(X+\nu)
 -
 \Vert
 X
 \Vert_{L^2(\Omega )}^2
 \Vert_{L^2(\Omega )}
 +
 \Vert
 \langle
 {D}_\cdot X , \EE \big[ D_\cdot X \mid \tilde{\cal F}_\cdot \big]
 \rangle
 -
 \EE[ \langle
 {D}_\cdot X, \EE \big[ D_\cdot X \mid \tilde{\cal F}_\cdot \big]
 \rangle
 ]
 \Vert_{L^2(\Omega )}
.
\end{eqnarray*}
\end{prop} 
\subsubsection*{Finite difference operator bound} 
Using the finite difference operator $\nabla$
we obtain the following bound which applies in particular
to multiple stochastic integrals, see Proposition~\ref{cor:dwI} below.
 \begin{prop}
   \label{thm:dWI}
 Let $X\in \Dom ( \nabla )$ be such that $\E[X]=0$. We have 
\begin{eqnarray} 
\label{dwf} 
\lefteqn{
  d_W(X,\mathcal N) \leq \E\[\left|1-\frac{1}{2} \langle\nabla_\cdot X, -\nabla_\cdot L^{-1}X\rangle\right|\]
}
\\
\nonumber
& &
 + \frac{1}{2}
    \E \left[ 
\int_0^\infty 
| \nabla_t L^{-1} X | 
| \nabla_t X |^2 dt
 \right] +  
\frac{1}{4}
    \E \left[ 
\int_0^\infty 
| \nabla_t L^{-1} X | 
\int_{2\lfloor t/2\rfloor}^{2\lfloor t/2\rfloor+2}
 | \nabla_s X |^2 ds dt
 \right]
.
\end{eqnarray} 
\end{prop}
\begin{Proof}
 By \eqref{eq:nablaaspsi}, for every function $f \in {\cal C}^2(\R )$, 
 the finite difference operator $\nabla$ satisfies 
\begin{eqnarray*} 
 \nabla_t f( X)
 & = & \frac{1}{2}
 \int_{2\lfloor t/2\rfloor}^{2\lfloor t/2\rfloor+2}
 ( f(X \circ \Psi_t) - f(X \circ \Psi_s) ) ds
 \\
 & = & 
 \frac{1}{2}
 \int_{2\lfloor t/2\rfloor}^{2\lfloor t/2\rfloor+2}
 \big(
 f'(X \circ \Psi_s) ( X \circ \Psi_t - X \circ \Psi_s )
 + R_f ( X \circ \Psi_t - X \circ \Psi_s )
 \big) ds
  \\
 & = & 
 \frac{1}{2}
 \int_{2\lfloor t/2\rfloor}^{2\lfloor t/2\rfloor+2}
 f'(X \circ \Psi_s) ( X \circ \Psi_t - X \circ \Psi_s )
 ds
 + \frac{1}{2}
 \int_{2\lfloor t/2\rfloor}^{2\lfloor t/2\rfloor+2}
 R_f ( X \circ \Psi_t - X \circ \Psi_s )
 ds, 
\end{eqnarray*} 
 $t \in \real_+$, where the function $R_f$ is
 such that $|R_f (y)|\leq y^2 \Vert f''\Vert_\infty /2$,
  $y\in \real$.
  Hence for any $f\in \mathcal{T}$,
  by the duality relation~\eqref{dr} we have
\begin{align}
  \nonumber
  & \E [ f'(X) - X f (X) ] 
=
\E [ f'(X) - X L L^{-1} f(X) ] 
  \\
\nonumber
  &=
  \E \left[ f'(X) - \frac{1}{2} \langle \nabla f (X) , - \nabla L^{-1} X\rangle \right] 
  \\
\nonumber
  &=
  \E \left[ f'(X) - \frac{1}{2} \int_0^\infty
    \nabla_t f (X) ( - \nabla_t L^{-1} X ) dt
    \right] 
  \\
\nonumber
  & = 
  \E \left[ f'(X) - \frac{1}{4} \int_0^\infty
 \int_{2\lfloor t/2\rfloor}^{2\lfloor t/2\rfloor+2}
 f'(X \circ \Psi_s) ( X \circ \Psi_t - X \circ \Psi_s )
 ds
 ( - \nabla_t L^{-1} X ) dt
    \right] 
  \\
\label{hfhjkf} 
  & - 
  \frac{1}{4} 
  \E \left[ \int_0^\infty
      \int_{2\lfloor t/2\rfloor}^{2\lfloor t/2\rfloor+2}
 R_f ( X \circ \Psi_t - X \circ \Psi_s )
 ds
  ( - \nabla_t L^{-1} X ) dt
    \right]. 
\end{align}
Regarding the first term, we note that for any two
square-integrable random variables $F$ and $G$,
by \eqref{eq:nablaaspsi} we have 
\begin{equation} 
  \label{jkdf}
  E\left[
 ( F \circ \Psi_t ) G 
 \right] 
=
\frac{1}{2} E\left[
 ( F \circ \Psi_t ) 
 \int_{2k}^{2k+2} 
 G \circ \Psi_s 
 ds
\right] 
\mbox{ and }
E\left[
  ( \nabla_t F )
  G 
 \right] 
=
\frac{1}{2} E\left[
  \nabla_t F 
  \int_{2k}^{2k+2} 
 G \circ \Psi_s 
 ds
\right] 
,
\end{equation} 
$t\in [2k,2k+2]$, $k\in \inte$, hence 
\begin{eqnarray*}
  \lefteqn{
    \! \! \!  \! \! \! \! \! \!  \! \! \!  \! \! \!  \! \! \!  \! \! \! 
    \left|
  \E \left[ f'(X) - \frac{1}{4} \int_0^\infty
 \int_{2\lfloor t/2\rfloor}^{2\lfloor t/2\rfloor+2}
 f'(X \circ \Psi_s) ( X \circ \Psi_t - X \circ \Psi_s )
 ds
 ( - \nabla_t L^{-1} X ) dt
    \right] 
  \right|
  }
  \\
  & = & 
  \left|
\E \left[ f'(X) - \frac{1}{2} \int_0^\infty
 f'(X ) ( X \circ \Psi_t - X )
 ( - \nabla_t L^{-1} X ) dt
    \right] 
    \right|
\\
  & = & 
  \left|
\E \left[ f'(X) \left( 1 - \frac{1}{2} \int_0^\infty
 ( X \circ \Psi_t - X )
 ( - \nabla_t L^{-1} X ) dt
  \right)
  \right] 
    \right|
  \\
  & = & 
  \left|
\E \left[ f'(X) \left( 1 - \frac{1}{2} \int_0^\infty
 \nabla_t X 
 ( - \nabla_t L^{-1} X ) dt
  \right)
  \right] 
    \right|
  \\
  & \leq & 
  \E \left[ \left|
    1 - \frac{1}{2} \int_0^\infty
 \nabla_t X 
 ( - \nabla_t L^{-1} X ) dt
  \right|
  \right] 
,
\end{eqnarray*}
because $\|f'\|_\infty \leq 1$.
Next, given that $\|f''\|_\infty \leq 2$, the term \eqref{hfhjkf} can be bounded as
\begin{eqnarray*}
  \lefteqn{
    \! \! \! \! \! \! \! \! \! \! \! 
    \frac{1}{4}
    \left|
  \E \left[ \int_0^\infty
      \int_{2\lfloor t/2\rfloor}^{2\lfloor t/2\rfloor+2}
 R_f ( X \circ \Psi_t - X \circ \Psi_s )
 ds
  ( - \nabla_t L^{-1} X ) dt
    \right] \right|
  }
  \\
   & \leq & 
\frac{1}{4}
    \E \left[ 
\int_0^\infty 
| \nabla_t L^{-1} X | 
\int_{2\lfloor t/2\rfloor}^{2\lfloor t/2\rfloor+2}
 | X \circ \Psi_t - X \circ \Psi_s |^2 
 ds dt
 \right]
  \\
   & = & 
\frac{1}{4}
    \E \left[ 
\int_0^\infty 
| \nabla_t L^{-1} X | 
\int_{2\lfloor t/2\rfloor}^{2\lfloor t/2\rfloor+2}
 | \nabla_t X - \nabla_s X |^2 
 ds dt
 \right]
  \\
   & = & 
\frac{1}{4}
    \E \left[ 
\int_0^\infty 
| \nabla_t L^{-1} X | 
\int_{2\lfloor t/2\rfloor}^{2\lfloor t/2\rfloor+2}
 ( | \nabla_t X |^2 + | \nabla_s X |^2 - 2 \nabla_s X \nabla_t X ) 
 ds dt
 \right]
  \\
   & = & 
\frac{1}{4}
    \E \left[ 
\int_0^\infty 
| \nabla_t L^{-1} X | 
\int_{2\lfloor t/2\rfloor}^{2\lfloor t/2\rfloor+2}
 ( | \nabla_t X |^2 + | \nabla_s X |^2 ) 
 ds dt
 \right]
  \\
   & = & 
 \frac{1}{2}
    \E \left[ 
\int_0^\infty 
| \nabla_t L^{-1} X | 
| \nabla_t X |^2 dt
 \right] +  
\frac{1}{4}
    \E \left[ 
\int_0^\infty 
| \nabla_t L^{-1} X | 
\int_{2\lfloor t/2\rfloor}^{2\lfloor t/2\rfloor+2}
 | \nabla_s X |^2 ds dt
 \right]
,
\end{eqnarray*} 
 where we used the relation 
$$
 E\left[
 ( F \circ \Psi_t ) 
 \int_{2k}^{2k+2} 
 \nabla_s G
 ds
\right] 
 =0
 \quad
 \mbox{and}
 \quad 
 E\left[
 \nabla_t F 
 \int_{2k}^{2k+2} 
 \nabla_s G
 ds
\right] 
 =0
, 
$$
 $t\in [2k,2k+2]$, $k\in \inte$,
 that hold similarly to \eqref{jkdf}. 
 We conclude to \eqref{dwf} by the inequality \eqref{fjhdfs},
 which is the bound $(2.33)$ in \cite{utzet2}. 
\end{Proof}
 The second term in \eqref{dwf} can also be written as 
\begin{eqnarray*}
  \lefteqn{
    \! \! \! \! \! \! \! \! \! \! \! \! \! \! \! \! \!
    \frac{1}{4}
\E \left[ 
\int_0^\infty 
| \nabla_t L^{-1} X | 
\int_{2\lfloor t/2\rfloor}^{2\lfloor t/2\rfloor+2}
 | X \circ \Psi_t - X \circ \Psi_s |^2 
 ds
dt
\right] 
  }
  \\
  & = &  
\frac{1}{4}
\E \left[ 
\int_0^\infty 
| \nabla_t L^{-1} X | 
\int_{2\lfloor t/2\rfloor}^{2\lfloor t/2\rfloor+2}
 | X \circ \Psi_t - X |^2 \circ \Psi_s  
 ds
dt
\right] 
  \\
  & = &  
\frac{1}{4}
\E \left[ 
\int_0^\infty 
| \nabla_t L^{-1} X | 
\int_{2\lfloor t/2\rfloor}^{2\lfloor t/2\rfloor+2}
 | X \circ \Psi_t - X |^2 ds
dt
\right] 
  \\
  & = &  
\frac{1}{2}
\E \left[ 
\int_0^\infty 
| \nabla_t L^{-1} X | 
 | X \circ \Psi_t - X |^2 
dt
\right] 
.
\end{eqnarray*}
 Taking $X= I_n(f_n)$ in Proposition~\ref{thm:dWI}, we get
the following result. 
\begin{prop}
  \label{cor:dwI}
  Let $f_n\in \hat{L}^2(\R_+^n)$. The following estimate holds: 
\begin{eqnarray*} 
  \lefteqn{
    \! \! \! 
    d_W( I_n(f_n),\mathcal N)
  \leq \sqrt{\E\[\(1-\frac1n\|\nabla I_n(f_n)\rangle\|^2_{L^2(\R_+,dx/2)}\)^2\]}
  }
  \\
\nonumber
& & +
\frac{1}{2n}
\E \left[ 
\int_0^\infty 
| \nabla_t I_n(f_n) |^3 
dt
\right]
+ \frac{1}{4n}
\E \left[ 
\int_0^\infty 
| \nabla_t I_n(f_n) | 
\int_{2\lfloor t/2\rfloor}^{2\lfloor t/2\rfloor+2}
 | \nabla_s I_n(f_n) |^2 
 ds
dt
\right]. 
\end{eqnarray*} 
\end{prop} 
\section{Single stochastic integrals}
\label{s2.1}
For single stochastic integrals, 
Proposition~\ref{cor:dwI} shows the following.
\begin{prop}
  \label{djkld}
  For $f_1\in L^2(\R_+)$ such that 
$\int_{2k}^{2k+2}f_1 (t)dt=0$, $k \in \inte$, 
    we have 
\begin{eqnarray} 
\nonumber 
  d_W( I_1(f_1),\mathcal N)
  & \leq & \left|1-\frac{1}{2} \int_0^\infty | f_1(t) |^2dt \right|
+
\frac{1}{2}
\int_0^\infty 
| f_1(t) |^3
dt
+
\frac{1}{4}
\int_0^\infty 
| f_1(t) | 
\int_{2\lfloor t/2\rfloor}^{2\lfloor t/2\rfloor+2}
| f_1(s) |^2 dsdt
\\
\label{eq:I_1nabla}
& \leq & \left|1-\frac{1}{2} \int_0^\infty | f_1(t) |^2dt \right|
+
\int_0^\infty 
| f_1(t) |^3
dt
.
\end{eqnarray} 
\end{prop}
 Consider now a sum $Z_n$ 
 of independent centered random variables $(X_k)_{k\geq 1}$ written,
 as in \eqref{zn}, as 
 \begin{equation}
   \label{zn2} 
Z_n : = I_1 \big( f_1 {\bf 1}_{[0,2n]} \big)
= \sum_{k=1}^n X_k, 
\end{equation} 
with $f_1 \in C^1_{\Box}(\R_+)\cap L^2(\R_+)$ given by \eqref{djkldd}
from the respective cumulative distribution functions $(F_k)_{k\geq 1}$.
In this case, Proposition~\ref{djkld} can be rewritten as follows. 
\begin{prop}
  \label{djkld2}
  Given $(Z_n)_{n \geq 1}$ written as in \eqref{zn2}
  we have 
  \begin{equation}
    \label{dww} 
    d_W( Z_n ,\mathcal N) 
  \leq 
 |1-E[Z_n^2]| +
 \sum_{k=1}^n E[|X_k|^3]
  +
\sum_{k=0}^\infty
E[|X_k|]E[|X_k|^2]
.
\end{equation} 
\end{prop} 
\begin{Proof} 
We note that $f_1(2k+1+U_k) = F^{-1}_k ((U_k+1)/2)$
has same distribution as $X_k$,
$k\geq 1$, hence \eqref{eq:I_1nabla} can be rewritten as 
\begin{eqnarray*} 
  \lefteqn{
    \! \! \! \! \! \! 
    d_W( Z_n ,\mathcal N) \leq |1-E[Z_n^2]| 
+
\frac{1}{2}
 \int_0^{2n} 
| f_1(t) |^3
dt
+
\frac{1}{4}
\sum_{k=0}^\infty
  \int_{2k}^{2k+2}
| f_1(t) | 
dt
\int_{2k}^{2k+2}
| f_1(s) |^2 ds 
  }
  \\
 & = & 
 |1-E[Z_n^2]| +
 \frac{1}{2}
 \sum_{k=1}^n \int_0^2
 | F^{-1}_k (t/2) |^3 dt
+
\frac{1}{4}
\sum_{k=0}^\infty
  \int_{0}^{2}
| F^{-1}_k (t/2) |
dt
\int_{0}^{2}
| F^{-1}_k (t/2) |^2 ds 
 \\
 & = & 
 |1-E[Z_n^2]| +
 \sum_{k=1}^n
 \int_{-\infty}^\infty 
 |x|^3 d F_k (x )
 +
\sum_{k=0}^\infty
 \int_{-\infty}^\infty 
 |x| d F_k (x )
 \int_{-\infty}^\infty 
 |y|^2 d F_k (y )
\\
  & = & 
 |1-E[Z_n^2]| +
 \sum_{k=1}^n E[|X_k|^3]
  +
\sum_{k=0}^\infty
E[|X_k|]E[|X_k|^2]
.
\end{eqnarray*} 
\end{Proof} 
Using H\"older's inequality,
Proposition~\ref{djkld2} shows that
\begin{equation}
   \label{ldjklf}
    d_W( \tilde{Z}_n ,\mathcal N)  \leq 
 \frac{2}{(E[(Z_n)^2])^{3/2}} 
 \sum_{k=1}^n E[|X_k|^3], \qquad n \geq 1, 
\end{equation} 
for the normalized sum
$\displaystyle 
\tilde{Z}_n : = 
 ( E[(Z_n)^2] )^{-1/2} \sum_{k=1}^n X_k
$, 
which recovers the bound \eqref{dw} of \cite{goldstein},
with however a worse constant. 
\subsubsection*{Bernoulli random variables} 
\begin{figure}[ht!]
 \centering
  \includegraphics[width=0.7\textwidth,height=0.3\textwidth]{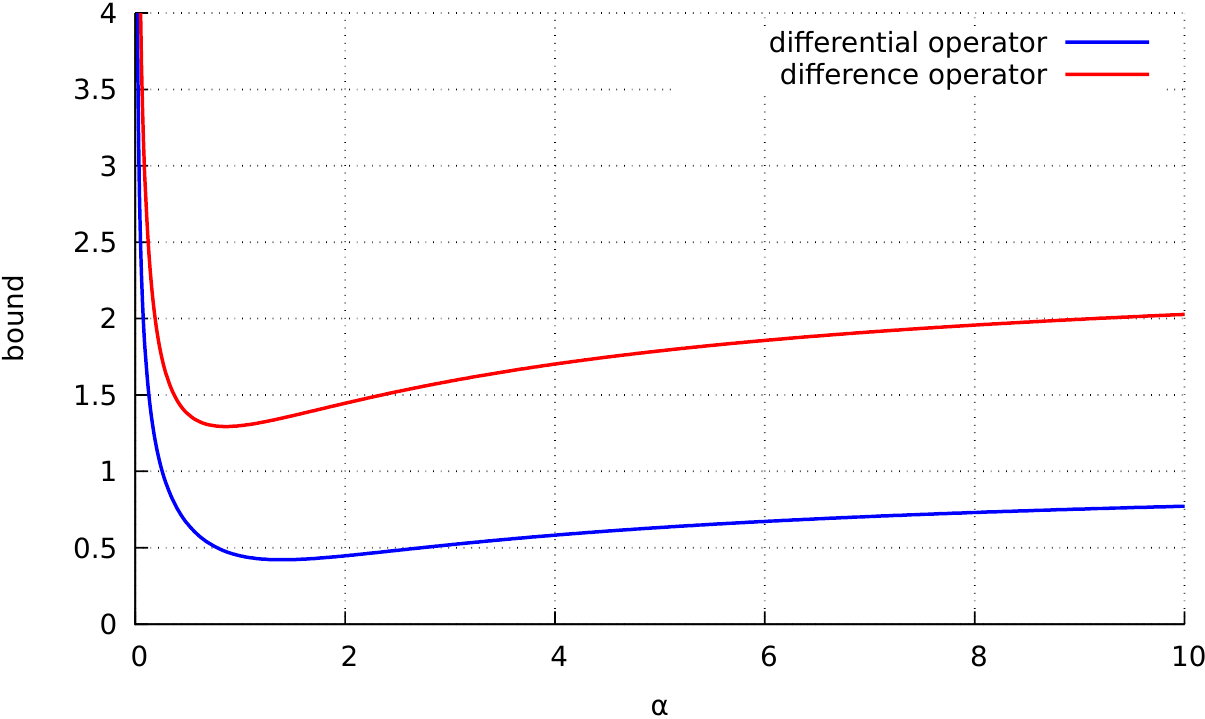}
   \vskip-0.4cm
   \caption{\small Comparison of bounds.}
  \label{f1}
\end{figure}

\begin{description}
  \item Given $(p_k)_{k\geq 1}$ a sequence in
$(0,1)$, letting 
$$
f(t) : = \sum_{k=1}^\infty
\frac{\alpha_k}{\sqrt{p_k(1-p_k)}}
\left( {\bf 1}_{[2k-2,2k-2+2p_k]}(t) - p_k \right),
\qquad t\in \real_+, 
$$
the single integral $I_1(f_1{\bf 1}_{[0,2n]})$ 
becomes a weighted sum
$$
I_1(f_1{\bf 1}_{[0,2n]})
 = \sum_{k=1}^n \alpha_k X_k
 $$ 
of centered and normalized Bernoulli random variables
$(X_k)_{k\geq 1}$ with parameters $(p_k)_{k\geq 1}$, and \eqref{dww} 
shows that
$$
d_W( I_1(f_1),\mathcal N) \leq
\left|1-\sum_{k=1}^\infty |\alpha_k|^2 
\right|
+
2 \sum_{k=1}^\infty | \alpha_k|^3 \frac{1 -2 p_k(1-p_k)}{\sqrt{p_k(1-p_k)}}, 
$$
which provides a simple distance bound for the
sum of non-symmetric Bernoulli random variables, cf. 
Corollary~3.3 of \cite{nourdin3}, 
Corollary~4.1 of \cite{privaulttorrisi4} and 
Theorem~4.1 of \cite{reichenbachs} for other versions. 
\end{description}
\noindent
By Proposition~\ref{prop:dw+dtv} we have the following
result that uses the derivation operator $D$.
\begin{prop}
  \label{djkdddfds}
  For $f_1 \in C^1_{\Box}(\R_+)\cap L^2(\R_+)$
  such that $\int_{2k}^{2k+2}f_1 (t)dt=0$, $k \in \inte$, 
    we have 
 \begin{eqnarray}
 \label{eq:I_1D}
  \lefteqn{
    d_W ( I_1(f_1),\mathcal N)
  }
  \\
  \notag
  & \leq&\left|1-\frac{1}{2}
  \int_0^\infty | f_1(t) |^2dt\right|
  +\frac12\sqrt{
    2 \int_0^\infty \left| f_1'(x) \int_0^{x}f_1(t)dt \right|^2dx-\sum_{k=0}^\infty\(\int_{2k}^{2k+2} | f_1(t) |^2dt\)^2}.  
\end{eqnarray} 
The bound for $d_{TV} ( I_1(f_1),\mathcal N)$ is twice as large
as \eqref{eq:I_1D}. 
\end{prop} 
\begin{Proof}
   We note that by \eqref{eq:uv} and Lemma~\ref{alm} we have 
\begin{eqnarray*}
  \lefteqn{
    \! \! \! \! \! \! \! \! \!
    \! \! \! \! \! \! \! \! \!
    E\left[
   \left(
    \langle
 D_\cdot Z_n
 ,
 \EE\big[ D_\cdot Z_n \mid \tilde{\cal F}_\cdot \big]
 \rangle
 -
 \EE[ \langle
 D_\cdot Z_n
 ,
 \EE\big[ D_\cdot Z_n \mid \tilde{\cal F}_\cdot \big]
 \rangle
 \right)^2 \right]
    }
    \\
    & = &
    \frac{1}{2} \int_0^{2n} \left| f'_1(x) \int_0^{x}f_1(t)dt \right|^2dx
    - \frac{1}{4} \sum_{k=1}^n \(\int_{2k-2}^{2k} | f_1(t) |^2dt\)^2, 
\end{eqnarray*}
  \end{Proof}
Proposition~\ref{djkdddfds} can be rewritten as follows
using sums $Z_n$ of random variables $(X_k)_{k\geq 1}$. 
  \begin{prop}
  \label{djkldddas}
   Assume that $(X_k)_{k\geq 1}$ is a sequence
  of independent centered random variables
  having non-vanishing continuous densities. 
  Then the sum 
$$
Z_n : = \sum_{k=1}^n X_k, \qquad n \geq 1, 
$$ 
 satisfies the bound 
\begin{equation}
   \label{dwdw} 
  d_W ( Z_n ,\mathcal N)
  \leq 
 |1 - E[Z_n^2]| 
 +
 \sqrt{ \sum_{k=1}^n
   \big( E [ ( \varphi_{X_k} (X_k) )^2 ] - (E[(X_k)^2])^2 \big)
  }
. 
\end{equation} 
The bound for $d_{TV} ( Z_n ,\mathcal N)$ is twice as large
as \eqref{dwdw}. 
\end{prop} 
  \begin{Proof}
     By Lemma~\ref{alm}, \eqref{eq:uv} and
 \eqref{djkf}, see (3.17) in \cite{viens}, we have 
\begin{eqnarray*}
  \lefteqn{
    \frac{1}{2} \int_0^{2n} \left| f'_1(x) \int_0^{x}f_1(t)dt \right|^2dx
    - \frac{1}{4} \sum_{k=1}^n \(\int_{2k-2}^{2k} | f_1(t) |^2dt\)^2
  }
  \\
 \nonumber
 & = & 
    \sum_{k=1}^n
    \left(
    \int_0^1 \left|
    \frac{1}{F_K'(F_k^{-1} (x)) }
 \int_0^x F_k^{-1} (t) dt \right|^2dx
  - ( E[(X_k)^2])^2 \right) 
      \\
  \nonumber 
  & = &
  \sum_{k=1}^n
    \left(
     \int_{-\infty}^\infty \left|
    \frac{1}{F_k'(y) }
 \int_{-\infty}^y x dF_k (x) \right|^2dF_k(y)
 - ( E[(X_k)^2])^2 \right)
       \\
  \nonumber 
  & = &
  \sum_{k=1}^n
   \big( E [ ( \varphi_{X_k} (X_k) )^2 ] - (E[(X_k)^2])^2 \big)
  , 
\end{eqnarray*}
\end{Proof}
 Next, we consider some particular cases.
 \begin{description}
   \item {\bf Gaussian case.} 
     The Stein kernel
     of $X_k$ centered Gaussian
 is given by 
$$ 
  \varphi_{X_k} ( y ) = 
 - \frac{1}{F'_k (y)} \int_{-\infty}^y x dF_k (x)
 = 
 E[X_k^2], \qquad y \in \real, \quad k\geq 1, 
$$
 and the bound \eqref{dwdw}
 recovers $d_W ( Z_n ,\mathcal N)\leq |1-E[Z_n^2]|$
 as expected. 
\item {\bf Gamma case.}
\noindent 
The Stein kernel of $X_k$ a centered gamma random
variable is 
 $\varphi_{Z_n} ( y ) = y+ E[Z_n^2]$, $y \in \real$, 
 hence
 $$
 E [ ( \varphi_{Z_n} (Z_n))^2 ] = E [ (Z_n+E[Z_n^2])^2 ] = E[Z_n^2](1+E[Z_n^2]),
 $$ 
 and the bound \eqref{dwdw} shows that the sum $Z_n$ satisfies 
$$ 
  d_W ( Z_n ,\mathcal N)
  \leq 
 |1 - E[Z_n^2]| 
 +
 \sqrt{E[Z_n^2]},
 \qquad n \geq 1. 
 $$
 By the scaling relation 
 $$\varphi_{a Z_n}(y) = a^2 \varphi_{Z_n}(y/a)
 = ay + a^2 E[Z_n^2], \qquad y \in \real,
$$
 we find that the normalized sum 
 $\tilde{Z}_n : = Z_n / \sqrt{E[Z_n^2]}$
 satisfies 
$$ 
  d_W ( \tilde{Z}_n ,\mathcal N)
  \leq 
  \frac{1}{\sqrt{E[Z_n^2]}}, \qquad n \geq 1. 
$$ 
  In particular, 
  in the i.i.d. case we
  have
$$ 
  d_W ( \tilde{Z}_n ,\mathcal N)
  \leq 
  \frac{1}{\sqrt{nE[X_1^2]}}, \qquad n \geq 1, 
$$ 
  which systematically improves on
\eqref{ldjklf} and on the bound \eqref{dw} of \cite{goldstein},
{\it i.e.}
\begin{eqnarray*} 
  \lefteqn{
    d_W ( Z_n ,\mathcal N)
  \leq 
 n \frac{E[|X_1|^3]}{( E[Z_n^2] )^{3/2} } 
  }
  \\
 & = & 
\frac{2}{\sqrt{n E[X_1^2]}}
\left(
\frac{2 \Gamma ( 3 + E[X_1^2], E[X_1^2]) + 2 (E[X_1^2])^{2+E[X_1^2]} e^{-E[X_1^2]} (1 + E[X_1^2]) }{
\Gamma ( 3 + E[X_1^2] ) } -1 
\right), 
\end{eqnarray*} 
where $\Gamma ( 3 + s, s)$ denotes the upper incomplete gamma function.
Indeed, the ratio
$$
2 
\left(
\frac{2 \Gamma ( 3 + s, s) + 2 s^{2 + s} e^{-s} (1 + s) }{
\Gamma ( 3 + s ) } -1 
\right), 
$$
 of the two bounds tends to infinity as $s$ tends to infinity, and has
smallest value $2$ as $s$ tends to $0$. 
\item {\bf Beta case.}
\noindent
When $X_k$ has the centered
Beta$(\alpha, 1)$ distribution, $\alpha >0$, $k\geq 1$, we have
$$
 F(x) = \left( \frac{\alpha}{\alpha + 1 } + x 
 \right)^\alpha,
 \qquad
 x\in \left[ 
 - 
 \frac{\alpha}{\alpha + 1}
 ,
 \frac{1}{\alpha + 1}
 \right], 
$$
 and $E[X_k^2] = \alpha / ((\alpha +1)^2(\alpha +2))$,
 hence by \eqref{djkd} we have 
 \begin{eqnarray}
   \label{qn3} 
  E [ ( \varphi_{X_k} (X_k))^2 ]
 & = & 
 \frac{1}{(\alpha+1)^2}
 E\left[ \left( \frac{\alpha}{\alpha +1} + X_k \right)^2
 \left( \frac{1}{\alpha +1} - X_k \right)^2 \right]
 \\
 \nonumber
 & = &
 \frac{2\alpha }{
   (\alpha + 4)(\alpha + 3 ) ( \alpha + 2)(\alpha+1)^2}, 
\end{eqnarray} 
 and by Proposition~\ref{djkldddas}, the normalized sum 
$$
 \tilde{Z}_n : = \frac{1}{\sqrt{n E[X_1^2]}} \sum_{k=1}^n X_k,
 \qquad n \geq 1, 
$$
 satisfies 
$$ 
  d_W ( \tilde{Z}_n ,\mathcal N)
  \leq
  \frac{1}{\sqrt{n}} 
  \sqrt{
    \frac{4+\alpha ( \alpha^2 + \alpha -2 ) }{\alpha  (\alpha +3) (\alpha +4)}},$$ 
 which systematically improves on
 \eqref{ldjklf} and on the bound \eqref{dw} of \cite{goldstein},
 {\it i.e.} 
$$
d_W ( Z_n ,\mathcal N)
\leq \frac{1}{\sqrt{n}} E[|X_1|^3]
= 
 \frac{2}{\sqrt{n}} \sqrt{\frac{\alpha +2}{\alpha }}
 \left(
 \frac{6 \alpha ( \alpha / (\alpha +1))^{\alpha + 1} + 1 - \alpha}{
   \alpha +3}
 \right), 
 $$
 as can be checked from Figure~\ref{f1}.
\end{description}

\vskip-0.4cm 

\noindent 
For example in the uniform case with $\alpha = 1$
we have $X_k = U_k$, $k\in \inte$,
and $F(x) = (x+1)/2$, $x\in [-1,1]$,
and 
$$f_1(t)= \sqrt{3}
\sum_{k=0}^\infty (t-2k-1)\mathbf1_{[2k,2k+2)}(t),
  $$ 
 hence \eqref{eq:I_1nabla} shows that 
 the sequence 
  $\displaystyle Z_n:=\sqrt{3/n} \sum_{k=1}^{n}U_k$,
     satisfies
  $$
  d_W(Z_n,\mathcal{N} )
  \leq
  \frac{3^{3/2}}{\sqrt{n}} E[|X_1|^3] 
  = 
  \frac{3}{4} \sqrt{\frac{3}{n}}
, 
$$
  whereas \eqref{eq:I_1D} yields $d_W (Z_n,\mathcal{N} ) \leq 
   1/\sqrt{5n}$. 
\section{Multiple stochastic integrals} 
\label{s3}
In this section we apply the multiplication formula given in the
appendix Section~\ref{s5} 
in order to obtain bounds on the
distance between multiple stochastic integrals and
the normal distribution $\mathcal N$.
 In the sequel for $0\leq i\leq k\leq n\wedge m$ 
  we define 
\begin{align}
  \label{star}
  & f_n\star_{k}^{i}g_m(\gamma_1,\ldots ,\gamma_{k-i}, t_1,\ldots ,t_{n-k},s_1,\ldots ,s_{m-k})\\
  \nonumber
  &: = 
\frac{1}{2^i}
\int_{[0,\infty)^i}\hskip-0.4cm
f_n(z_1,\ldots ,z_i,\gamma_1,\ldots ,\gamma_{k-i}, t_1,\ldots ,t_{n-k}) g_m(z_1,\ldots ,z_i,\gamma_1,\ldots ,\gamma_{k-i},s_1,\ldots ,s_{m-k})dz_1\cdots dz_i, 
\end{align}
 and we denote by $f_n \hskip0.1cm \widetilde{\star}_{k}^{i} g_m$ the symmetrisation of $f_n\star_{k}^{i}g_m$, i.e.
$$ 
f_n \hskip0.1cm \widetilde{\star}_{k}^{i} g_m(x_1,\ldots ,x_{m+n-k-i})
:=\frac{1}{(m+n-k-i)!}\sum_{\sigma\in S_{m+n-k-i}}f_n\star_{k}^{i}g_m(x_{\sigma(1)},\ldots ,x_{\sigma(m+n-k-i)}),
$$ 
where $S_{m+n-k-i}$ stands for the group of all permutations of the set $\{1,\ldots ,m+n-k-i\}$.
Note that $f_n\star_{k}^{i}g_m$ may not satisfy \eqref{ass:int0}, 
even if $f_n$ and $g_m$ satisfy \eqref{ass:int0}. 
The multiplication formula of Theorem~\ref{thmm} below 
can be given in many different forms,
 one of which is presented in the next Proposition~\ref{cor:mf}.
\begin{prop}
\label{cor:mf}
Let  $f_n$, $g_m$ satisfy \eqref{ass:int0} and $f_n\star_{k}^{i}g_m\in L^2(\R_+^{m+n-k-i})$ for every $0\leq i \leq k\leq m\wedge n$.
We have
$$ I_n(f_n) I_m(g_m)=\sum_{k=0}^{2(m\wedge n)} I_k(h_k),$$
where 
$$h_k=\sum_{r=0}^{n\wedge m}\sum_{l=0}^r\mathbf1_{\{ 2n-r-l=k \}}\,r!{{n}\choose{r}}{{m}\choose{r}}{{r}\choose{l}}f_n\hskip0.1cm \widetilde{\star}^l_r g_m.$$
\end{prop}
\subsubsection*{Bounds obtained from the finite difference operator $\nabla$}
To obtain a more explicit bound than
in Proposition~\ref{cor:dwI}
 we have to employ the multiplication formula.
Precisely, by virtue of
Proposition~\ref{cor:mf} we may express $\big(I_n(f_n)\big)^2$ as follows:
\begin{equation}\label{eq:G}
\big( I_n(f_n)\big)^2=\sum_{k=0}^{2n} I_k\big(G^n_kf_n\big),
\end{equation}
where
$$G^n_kf_n(z_1,\ldots ,z_k)=
{\bf 1}_{\Delta_k} (z_1,\ldots , z_k)
\sum_{r=0}^{n}\sum_{l=0}^r\mathbf1_{\{ 2n-r-l=k \}}r!{{n}\choose{r}}^2{{r}\choose{l}}f_n\hskip0.1cm \widetilde{\star}^l_rf_n (z_1,\ldots ,z_k).$$
\begin{corollary}
  \label{thm:dIZnabla}
  Let $f_n\in L^2(\R_+^n)$ be a symmetric function
  satisfying \eqref{ass:int0}. 
  Assume that 
$$\widehat G^{n}_kf_n (*) : =
\frac{1}{2} \int_0^\infty G_k^{n-1}f_n(t,*)dt
$$
belongs to $\in L^2(\R_+^k)$ for all
$1\leq  k\leq 2n-2$. We have
\begin{eqnarray*}
  \lefteqn{ 
    \! \! \! \! \! \! \! \! \! \! \! \! 
    d_W( I_n(f),\mathcal N)\leq \sqrt{\( 1 - n!\left\|f_n\right\|^2_{L^2(\R_+,dx/2)^{\circ n}} \)^2
  +n^2\sum_{k=1}^{2n-2}k! \left\| \widehat G^{n}_kf_n\right\|^2_{L^2(\R_+,dx/2)^{\circ k}}}
  }
  \\ 
  & & + n^2\sqrt{2(n-1)!} \left\|f_n \right\|_{L^2(\R_+,dx/2)^{\circ n}} 
  \sqrt{\sum_{k=0}^{2n-2}k!\int_0^\infty\left\|G^{n-1}_kf_n(t,\cdot)\right\|^2_{L^2(\R_+,dx/2)^{\circ k}}dt}.
\end{eqnarray*}
\end{corollary} 
\begin{Proof}
  We are going to estimate both components appearing in
  Proposition~\ref{cor:dwI}. 
The formula \eqref{eq:G} lets us write 
\begin{align*}
\big(\nabla_t I_n(f_n)\big)^2 =nn!\left\|f_n(t,\cdot)\right\|^2_{L^2(\R_+,dx/2)^{\circ (n-1)}}+n^2\sum_{k=1}^{2n-2} I_k\big(G^{n-1}_kf_n(t,\cdot)\big).
\end{align*}
Hence we have 
$$ 
\frac1n \|\nabla_\cdot I_n(f_n) \|^2_{L^2(\R_+,dx/2)}
=n!\left\|f_n\right\|^2_{L^2(\R_+,dx/2)^{\circ n}}
+\frac{n}{2}
\sum_{k=1}^{2n-2}\int_0^\infty I_k\big(G^{n-1}_kf_n(t,\cdot)\big)dt.
$$ 
Since multiple integrals of different orders are orthogonal, we get
\begin{eqnarray*}
  \lefteqn{ 
    \E\[\(1-\frac1n\|\nabla_\cdot I_n(f_n)\rangle\|^2_{L^2(\R_+,dx/2)}\)^2\]
  }
  \\
   & = & \( 1 - n!\left\|f_n\right\|^2_{L^2(\R_+,dx/2)^{\circ n}} \)^2
+n^2\sum_{k=1}^{2n-2}\E\left[\(\int_0^\infty I_k\big(G^{n-1}_kf_n(t,\cdot)\big)dt/2\)^2\right].
\end{eqnarray*}
Finally, by \eqref{eq:EI^2}, we obtain
\begin{align*}
\E\left[\left| \int_0^\infty I_k\big(G^{n-1}_kf_n(t,\cdot)\big)dt/2\right|^2\right]=\E\left[\( I_k\big(\widehat G^{n}_kf_n\big)\)^2\right]\leq k! \left\| \widehat G^{n}_kf_n\right\|^2_{L^2(\R_+,dx/2)^{\circ k}},
\end{align*}
 which implies
$$ 
    \E\[\left| 1-\frac1n\|\nabla_\cdot I_n(f_n)\rangle\|^2_{L^2(\R_+,dx/2)}\right|^2\]
 \leq
  \left| 1 - n!\left\|f_n\right\|^2_{L^2(\R_+,dx/2)^{\circ n}} \right|^2
+n^2\sum_{k=1}^{2n-2}k! \left\| \widehat G^{n}_kf_n\right\|^2_{L^2(\R_+,dx/2)^{\circ k}}.
$$ 
To get the second component of the estimates in the thesis we use Cauchy-Schwarz inequality in the following way: 
\begin{eqnarray*} 
  \lefteqn{
    \! \! \! \! \! \! \! \! \! \! \! \! \! \! \! \! \! \! \! \! \! \! \! \! \! 
    \frac1{2n} \int_0^\infty\E\[|\nabla_t I_n(f_n)|^3\]dt
    \leq \frac{1}{2}
    \sqrt{\int_0^\infty\E\[\big( I_{n-1}(f_n(t,*))\big)^2\]dt}\sqrt{\int_0^\infty\E\[|\nabla_t I_n(f_n)|^4\]dt}
  }
  \\
  &\leq & n^2
  \left\|f_n \right\|_{L^2(\R_+,dx/2)^{\circ n}}
  \sqrt{\frac{(n-1)!}{2} \int_0^\infty\E\[\big( I_{n-1}(f_n(t,*))\big)^4\]dt }.
\end{eqnarray*}
Since 
\begin{align*}
\big( I_{n-1}(f_n(t,*))\big)^2=\sum_{k=0}^{2n-2} I_k\big(G^{n-1}_kf_n(t,\cdot)\big),
\end{align*}
and by orthogonality of multiple integrals, we have
\begin{align*}
\int_0^\infty\E\big[\big( I_{n-1}(f_n(t,*))\big)^4\big]dt\leq \sum_{k=0}^{2n-2}k!\int_0^\infty\left\|G^{n-1}_kf_n(t,\cdot)\right\|^2_{L^2(\R_+,dx/2)^{\circ k}}dt,
\end{align*}
which ends the proof.
\end{Proof}
As noted above, $I_n(f_n)$ can be used to represent
various U-statistics, including polynomials of Bernoulli
random variables, in which case Corollary~\ref{thm:dIZnabla}
provides an alternative to the results of
\cite{nourdin3}, \cite{krokowski}, 
\cite{privaulttorrisi4} for Bernoulli processes.
\subsubsection*{Bounds obtained from the derivation operator $D$}
Here we let $C^1_{\Box}(\R_+^n)$ denote the set of
functions which are ${\cal C}^1$ on every set of the form
$$
(2k_1,2k_1+2)\times\cdots \times(2k_n,2k_n+2),
\quad
k_1,\ldots ,k_n\in\N.
$$ 
Given $f_n \in \mathcal C^1_{\Box}(\R_+^n) \cap L^2(\R_+^n)$,
we define 
\begin{align*}
  & \! \! \! \! \! \! 
  H_k(s,z_2,\ldots ,z_{k+1}) := 
  \\
  & 
  \sum_{r=0}^{ n-1}\sum_{l=0}^r\mathbf1_{\{ 2n-2-r-l=k \}}\,r!{{n-1}\choose{r}}^2{{r}\choose{l}}
 \(    \partial_1f_n(s,*)\hskip0.1cm \widetilde{\star}^l_r \int_0^s f_n(t,*)\mathbf1_{\{*<s \}}dt\) 
(z_2,\ldots ,z_{k+1}),
\end{align*} 
 and
\begin{align*}
  & \! \! \! \! \! \! \! \! \! \! \! \! 
  J_k(s,z_1,\ldots ,z_k)
  \\
   & :=\sum_{r=0}^{ n-1}\sum_{l=0}^r\mathbf1_{\{ 2n-2-r-l=k \}}\,r!{{n-1}\choose{r}}^2{{r}\choose{l}}
  \( f_n(s,*)\widetilde{\star}^l_rf_n(s,*) \mathbf1_{\{*<s\}}\) (z_1,\ldots ,z_{k}),
\end{align*}  
 $1\leq k\leq 2n-2$, where
  $$
   \{*<u\} := \{x\in\R^{n-1}\ : \ x_i<u, \ i=1,\ldots ,n-1\}, 
  $$
   and assume that $H_k, J_k\in L^2(\R_+^{k+1})$.
 Additionally, we denote
$$\widehat J_k(z_1,\ldots ,z_k)=\frac{1}{2}
\int_0^\infty J_k(s,z_1,\ldots ,z_k)ds.$$
Next is a consequence of Proposition~\ref{prop:dw+dtv}.
  \begin{corollary}\label{thm:dIZD}
  Let $f_n \in \mathcal C^1_{\Box}(\R_+^n) \cap L^2(\R_+^n)$
  and satisfy \eqref{ass:int0}. We have
\begin{align*}
& 
   d_W ( I_n(f),\mathcal N)
    \leq \left|1-n!\left\|f_n\right\|^2_{L^2(\R_+,dx/2)^{\circ n}}\right|
  \\
  &
+n^2
\sqrt{\frac{1}{2}
\sum_{k=0}^{2n-2}\int_0^\infty\hskip-0.2cm
\E\[| I_k(H_k( s ,*)) |^2\]ds
+\sum_{k=1}^{2n-2}\E\big[\big| I_k(\widehat J_k)\big|^2\big]
-\frac{1}{4} \sum_{i=0}^\infty\sum_{k=0}^{2n-2} \E\[\left| I_k\(\int_{2i}^{2i+2}\hskip-0.4cm
J_k(s,*)ds \)\right|^2\]}
\\[10pt]
&\leq \left|1-n!\left\|f_n\right\|^2_{L^2(\R_+,dx/2)^{\circ n}}\right|+n^2\Bigg(
  \frac{1}{2}
  \sum_{k=0}^{2n-2}k!\int_0^\infty
  \left\|H_k( s ,*)\right\|^2_{L^2(\R_+,dx/2)^{\circ k}}ds
 +\sum_{k=1}^{2n-2}k!\left\|\widehat J_k\right\|^2_{L^2(\R_+,dx/2)^{\circ k}}\\
& + \left( (n-1)!\right)^2\sum_{i=0}^\infty\left| \int_{2i}^{2i+2}\left\|f_n(s,*) \mathbf1_{\{*<s\}}\right\|^2_{L^2(\R_+,dx/2)^{\circ (n-1)}}ds\right|^2\Bigg)^{1/2}.
\end{align*}
The bounds for $d_{TV}( I_n(f),\mathcal N)$ are equal to those for $d_W( I_n(f),\mathcal N)$ multiplied by 2.
\end{corollary}
\begin{Proof}
 By Lemma \ref{lem:cov} and formula~\eqref{eq:EI^2} we get
$$ 
 \E \big[ \langle D_\cdot I_n(f_n) , \E [ D_\cdot I_n(f_n)\mid
 \tilde{{\cal F}_\cdot } ]\rangle \big]
 =2\E\[( I_n(f_n))^2\]=2 n!\left\|f_n\right\|^2_{L^2(\R_+,dx/2)^{\circ n}}.
$$ 
 Next, we are going to provide an explicit form for the expression
$\langle D_\cdot I_n(f_n) ,
\E\big[ D_\cdot I_n(f_n) \mid \tilde{{\cal F}_\cdot }\big]\rangle$. We have
\begin{align*}
& 
  D_t I_n(f_n)
  \\
  &=n\sum_{0 \leq k_1\neq\cdots \neq k_{n-1}}\left(
  (1-U_{\lfloor t/2 \rfloor})\mathbf1_{(2{\lfloor t/2 \rfloor},2{\lfloor t/2 \rfloor}+1+U_{\lfloor t/2 \rfloor}]}(t)-(1+U_{\lfloor t/2 \rfloor})\mathbf1_{(2{\lfloor t/2 \rfloor}+1+U_{\lfloor t/2 \rfloor},2{\lfloor t/2 \rfloor}+2]}(t)
      \right)
      \\
& \qquad \times\partial_1f_n(2{\lfloor  t/2\rfloor}+1+U_{{\lfloor t/2 \rfloor}},2k_1+1+U_{k_1},\ldots , 2k_{n-1}+1+U_{k_{n-1}})\\
      &= n\left(
      (1-U_{\lfloor t/2 \rfloor})\mathbf1_{(2{\lfloor t/2 \rfloor},2{\lfloor t/2 \rfloor}+1+U_{\lfloor t/2 \rfloor}]}(t)-(1+U_{\lfloor t/2 \rfloor})\mathbf1_{(2{\lfloor t/2 \rfloor}+1+U_{\lfloor t/2 \rfloor},2{\lfloor t/2 \rfloor}+2]}(t)\right)
          \\
& \qquad \times I_{n-1}\(\partial_1f_n(2{\lfloor  t/2\rfloor}+1+U_{{\lfloor t/2 \rfloor}},*)\).
\end{align*}
By Proposition~10 and Lemma~1 in \cite{prebub} we get
\begin{align*}
\E\big[ D_t I_n(f_n) \mid \tilde{{\cal F}}_{t}\big]=\E\big[\nabla_t I_n(f_n) \mid \tilde{{\cal F}}_{t}\big]=n I_{n-1}\(f_n(t,*)\mathbf1_{\{*<2\lfloor t/2\rfloor\}}\).
\end{align*}
 Consequently, using the assumption \eqref{ass:int0} twice, we arrive at
\begin{align*} 
  &\langle D_\cdot I_n(f_n) ,
  \E\big[ D_\cdot I_n(f_n) \mid \tilde{{\cal F}_\cdot }\big] \rangle\\
  &=n^2
  \int_0^\infty \left(
  (1-U_{\lfloor t/2 \rfloor})\mathbf1_{(2{\lfloor t/2 \rfloor},2{\lfloor t/2 \rfloor}+1+U_{\lfloor t/2 \rfloor}]}(t)-(1+U_{\lfloor t/2 \rfloor})\mathbf1_{(2{\lfloor t/2 \rfloor}+1+U_{\lfloor t/2 \rfloor},2{\lfloor t/2 \rfloor}+2]}(t)
      \right) 
    \\
    & 
    \qquad \times I_{n-1}\(\partial_1f_n(2{\lfloor  t/2\rfloor}+1+U_{{\lfloor t/2 \rfloor}},*)\) I_{n-1}\(f_n(t,*)\mathbf1_{\{*<2\lfloor t/2\rfloor\}}\)dt
      \\
      &=n^2
      \int_0^\infty
      \left( \mathbf1_{(2{\lfloor t/2 \rfloor},2{\lfloor t/2 \rfloor}+1+U_{\lfloor t/2 \rfloor}]}(t)-\mathbf1_{(2{\lfloor t/2 \rfloor}+1+U_{\lfloor t/2 \rfloor},2{\lfloor t/2 \rfloor}+2]}(t)\right)
                    \\
          &
          \qquad \times I_{n-1}\(\partial_1f_n(2{\lfloor  t/2\rfloor}+1+U_{{\lfloor t/2 \rfloor}},*)\) I_{n-1}\(f_n(t,*)\mathbf1_{\{*<2\lfloor t/2\rfloor\}}\)dt
          \\
          &=2n^2\int_0^\infty
          \left(
          \mathbf1_{(2{\lfloor t/2 \rfloor},2{\lfloor t/2 \rfloor}+1+U_{\lfloor t/2 \rfloor}}(t) I_{n-1}\(\partial_1f_n(2{\lfloor  t/2\rfloor}+1+U_{{\lfloor t/2 \rfloor}},*)\) I_{n-1}\(f_n(t,*)\mathbf1_{\{*<2\lfloor t/2\rfloor\}}\)\right)
          dt
          \\
          &=2n^2\sum_{k=0}^\infty I_{n-1}\(\partial_1f_n(2k+1+U_{k},*)\) I_{n-1}\(\int_0^{2k+1+U_{k}}f_n(t,*)\mathbf1_{\{*<2k\}}dt\)
          \\
&=2n^2\int_0^\infty h(s)d(Y_s-s/2)+n^2\int_0^\infty h(s)ds,
\end{align*}
where
$$
h(s):= I_{n-1}\(\partial_1f_n(s,*)\) I_{n-1}\(\int_0^{s}f_n(t,*)\mathbf1_{\{*<s\}}dt\),
\qquad s\in\real_+, 
$$
is a random process when $n\geq 2$.
Note that by integration by parts we have 
\begin{align*}
\int_0^\infty h(s)ds=&\int_0^\infty I_{n-1}\(f_n(s,*)\) I_{n-1}\(f_n(s,*)\mathbf1_{\{*<s\}}\)ds.
\end{align*}
By the Fubini theorem we may express $\int_0^\infty h(s)d(Y_s-s/2)$ and $\int_0^\infty h(s)ds$ as $2n-1$ and $2n-2$ integrals with respect to $d(Y_s-s/2)$, respectively.  Then, applying \eqref{eq:uv} $2n-2$ times together with \eqref{ass:int0}, we obtain
\begin{align*}
& \E \[\int_0^\infty h(s)d(Y_s-s/2)\int_0^\infty h(s)ds\]
  \\
&= \int_{\R_+^{2n-1}}\E\[\int_0^\infty \partial_1f_n(u,x_1,\ldots ,x_{n-1})\int_0^uf_n(t,x_n,\ldots ,x_{2n-2})dt\mathbf1_{\{ (x_n,\ldots ,x_{2n-2})<u \}}d(Y_u-u/2)\]\\
&  \times
f_n(s,x_1,\ldots ,x_{n-1})f_n(s,x_n,\ldots ,x_{2n-2})\mathbf1_{\{ (x_n,\ldots ,x_{2n-2})<s\}}dx_1\cdots dx_{2n-2}ds
\\
 & =  0,
\end{align*}
and consequently
$$ 
\E \left[ \langle D_\cdot I_n(f_n) , \E\big[ D_\cdot I_n(f_n) \mid \tilde{{\cal F}_\cdot }\big]\rangle^2\right]
=4n^4\E\[\(\int_0^\infty h(s)d(Y_s-s/2)\)^2\]+n^4\E\[\(\int_0^\infty h(s)ds\)^2\].
$$ 
Using the orthogonality of multiple integrals of different orders
 and the relation 
$$ 
I_{n-1} ( f_n(s,*)) I_{n-1} ( f_n(s,*) ) \mathbf1_{\{*< s \}}=\sum_{k=0}^{2n-2} I_k\(J_k(s,*)\), 
$$ 
 we rewrite the latter component as follows: 
\begin{eqnarray*}
  \lefteqn{ 
  \! \! \! \! \! \! \!  \E\[\(\int_0^\infty h(s)ds\)^2\]=4\E\[\(\sum_{k=0}^{n-1} I_k(\widehat J_k)\)^2\]
  }
  \\
&=&4\sum_{k=1}^{n-1}\E\[\( I_k(\widehat J_k)\)^2\]+
\left( (n-1)!\int_0^\infty \left\|f_n(s,*)\mathbf1_{\{*<s\}}\right\|^2_{L^2(\R_+,dx/2)^{\circ (n-1)}}ds \right)^2\\
&=&4\sum_{k=1}^{n-1}\E\[\( I_k(\widehat J_k)\)^2\]+ \frac{4}{n^2}
 \left\|f_n\right\|^4 _{L^2(\R_+,dx/2)^{\circ n}}.
\end{eqnarray*}
Furthermore, by Proposition~\ref{cor:mf} we have
$$ 
  I_{n-1} ( \partial_1f_n( s ,*))
  I_{n-1} \left(
  \int_0^s f_n(t,*)\mathbf1_{\{*< s \}}dt
  \right)
  =\sum_{k=0}^{2n-2} I_k\(H_k( s ,*)\), 
$$
 hence \eqref{eq:uv} gives us
\begin{align*}
\E&\[\(\int_0^\infty h(s)d(Y_s-s/2)\)^2\]= \frac{1}{2}
\E\[\int_0^\infty |  h(s) |^2 ds\]-
\frac{1}{4}
\sum_{i=0}^\infty\E\[ \(\int_{2i}^{2i+2}h(s)ds \)^2\]\\
&=
\frac{1}{2}
\E\[\int_0^\infty\( \sum_{k=0}^{2n-2}I_k(H_k( s ,*))\)^2ds\]-\sum_{i=0}^\infty \E\[\(\sum_{k=0}^{n-1} I_k\(\int_{2i}^{2i+2}J_k(s,*)ds/2\)\)^2\]\\
&= \frac{1}{2}
\sum_{k=0}^{2n-2}\int_0^\infty\E\[\( I_k(H_k(s ,*))\)^2\]ds
-\sum_{i=0}^\infty\sum_{k=0}^{n-1} \E\[\( I_k\(\int_{2i}^{2i+2}J_k(s,*)ds/2\)\)^2\].
\end{align*}
We apply this to Proposition~\ref{prop:dw+dtv} and get the first inequality in the assertion of the  theorem. In order to derive the other one we use \eqref{eq:I^2<} and the estimate 
\begin{eqnarray*}
  \lefteqn{
    \! \! \! \! \! \! \! \! \! \! \! \! \! \! \! \! 
  \sum_{k=0}^{n-1} \E\[\( I_k\(\int_{2i}^{2i+2}J_k(s,*)ds/2\)\)^2\]
  \geq  \E\[\( I_0\(\int_{2i}^{2i+2}J_0(s,*)ds/2\)\)^2\]
  } 
  \\
  &= & ( (n-1)! )^2\(\int_{2i}^{2i+2}\left\|f(s,*)\mathbf1_{\{*<s\}}\right\|^2_{L^2(\R_+,dx/2)^{\circ (n-1)}}ds\)^2.
\end{eqnarray*}
\end{Proof}
\subsubsection*{A combinatorial central limit theorem}
In this section, we show that the
bounds of \cite{reichenbachs} for the 
Rademacher combinatorial central limit theorem of \cite{blei}
can be extended to our setting of random sequences.
\\
 
Given $K$ a symmetric subset of
$\tilde \Delta_q : =\{a\in\N^q:a_i\neq a_j\ if\ i\neq j\}$,
$(b_k)_{k\geq 0}$ a sequence
of real numbers, and $(X_k)_{k\geq 0}$ an i.i.d. sequence of
random variables such that $\E[X_1]=0$ and $\E[X_1^2]<\infty$,
define
$$S^{(b)}(K) : =\frac{1}{(q!\mu_b^{\otimes q}(K)(\E[X^2])^q)^{1/2}}\sum_{(i_1,\ldots ,i_q)\in K}b_{i_1}\cdots b_{i_q}X_{i_1}\cdots X_{i_q}.$$
Following \S~6.3 of \cite{reichenbachs}, 
we let $K^*_j$ denote the collection of all $(i_1,\ldots ,i_q) \in K$ such that $i_k = j$ for some $k \in\{1,\ldots ,q\}$,
and we define $K^\# \subset K\times K$
by stating that a pair $(i_1,\ldots ,i_q),(j_1,\ldots ,j_q)$ belongs to $K^\#$ if $\{i_1,\ldots ,i_q\}\cap\{j_1,\ldots ,j_q\} = \phi$ and there are $(k_1,\ldots ,k_q),(l_1,\ldots ,l_q) \in K$ such that $\{k_1,\ldots ,k_q,l_1,\ldots ,l_q\} = \{i_1,\ldots ,i_q,j_1,\ldots ,j_q\}$ and
 $(k_1,\ldots ,k_q)$ does not coincide with $(i_1,\ldots ,i_q)$ or $(j_1,\ldots ,j_q)$. 
\begin{theorem}
  \label{theoremk}
  There exists a constant $C=C(q)$ such that
$$d_W \(S^{(b)}(K),\mathcal N\)\leq C\(\E[X_1^4]\)^q\[\frac{\(\mu_b^{\otimes(2q)}(K^\#)\)^{1/2}}{\mu_b^{\otimes q}(K)}+\(\sup_{j\geq 1}\frac{\mu_b^{\otimes q}(K^*_j)}{\mu_b^{\otimes q}(K)}\)^{1/4}\].$$
\end{theorem}
\begin{Proof}
  Let $F$ be the distribution function of $X_1$
  with generalised inverse function $F^{-1}$. Then 
 we have $S^{(b)}(K)\stackrel{d}{=}I_q(f_q)$, 
where
$$f_q(t_1,\ldots ,t_q)=\frac{\mathbf1_{K}(\lfloor t_1/2\rfloor,\ldots ,\lfloor t_q/2\rfloor)}{(q!\mu_b^{\otimes q}(K)(\E[X_1^2])^q)^{1/2}}b_{\lfloor t_1/2\rfloor}\cdots b_{\lfloor t_q/2\rfloor}F^{-1}\(\frac{t_1-2\lfloor t_1/2\rfloor}2\)
\cdots F^{-1}\(\frac{t_q-2\lfloor t_q/2\rfloor}2\).$$
By Theorem  \ref{thm:dIZnabla}, there exist constants
$C_1$, $C_2$ depending only on $q$, such that
\begin{align}\nonumber
    d_W( I_q(f_q),\mathcal N)&\leq C_1\(
  \sum_{k=1}^{2q-2} \left\| \widehat G^{q}_kf_q\right\|_{L^2(\R_+,dx/2)^{\circ k}}+ 
 \sum_{k=0}^{2q-2}\int_0^\infty\left\|G^{q-1}_kf_q(t,\cdot)\right\|_{L^2(\R_+,dx/2)^{\circ k}}dt\)
\\\label{eq:aux3}
&\leq C_2\(
\sum_{k=1}^{q-1}\sum_{r=1}^k \left\| \(f_q\star_k^rf_q\)
    {\bf 1}_{\Delta_{2q-k-r}}
\right\|_{L^2(\R_+)^{\circ (2q-k-r)}}+ 
 \sum_{k=1}^q \sum_{r=0}^{k-1} \left\| f_q\star_k^rf_q\right\|_{L^2(\R_+)^{\circ (2q-k-r)}}\).
\end{align}
Note that for $r \leq k$ we have
\begin{align*}
  &\left\| \(f_q\star_k^rf_q\){\bf 1}_{\Delta_{2q-k-r}} \right\|_{L^2(\R_+)^{\circ (2q-k-r)}}
  =2^{q+r}\(\int_0^2\(F^{-1}(s)\)^2ds\)^{2r+2q-2k}\(\int_0^2\(F^{-1}(s)\)^4ds\)^{k-r}
  \\
  &
  \quad \times \sum_{y,z\in\N^{q-k}}\sum_{x\in\N^{k-r}}\(\sum_{\substack{(w_1,\ldots ,w_r,x_1,\ldots ,x_{k-r},y_1,\ldots ,y_{q-k})\in K\\(w_1,\ldots ,w_r,x_1,\ldots ,x_{k-r},z_1,\ldots ,z_{q-k})\in K}}\tilde f_q(w,x,y)\tilde f_q(w,x,z)\)^2{\bf 1}_{\tilde \Delta_{2q-k-r}}(x,y,z)
  \\
  &\leq2^{2q}\(\E[X_1^4]\)^q\left\| \(\tilde f_q\tilde\star_k^r\tilde f_q\){\bf 1}_{\tilde \Delta_{2q-k-r}}
  \right\|_{l^2(\N)^{\circ (2q-k-r)}},
\end{align*}
where, as in \cite{reichenbachs} or \cite{nourdin3},
the notation $\tilde \star$
is here the discrete version of the product
defined in \eqref{star}, and 
$$\tilde f_q(i_1,\ldots ,i_q) := \frac{\mathbf1_{K}(i_1,\ldots ,i_q)}{(q!\mu_b^{\otimes q}(K)(\E[X^2])^q)^{1/2}}b_{i_1}\cdots b_{i_q}
.$$
Furthermore, for fixed $y,z\in\N^{q-k}$ we get 
\begin{align*}
  &\left\| \(\tilde f_q\tilde\star_k^r\tilde f_q\){\bf 1}_{\tilde \Delta_{2q-k-r}}
  \right\|_{l^2(\N)^{\circ (2q-k-r)}}\\
  & \quad =\sum_{y,z\in\N^{q-k}}\sum_{x\in\N^{k-r}}\(\sum_{\substack{(w_1,\ldots ,w_r,x_1,\ldots ,x_{k-r},y_1,\ldots ,y_{q-k})\in K\\(w_1,\ldots ,w_r,x_1,\ldots ,x_{k-r},z_1,\ldots ,z_{q-k})\in K}}\tilde f_q(w,x,y)\tilde f_q(w,x,z)\)^2{\bf 1}_{\tilde \Delta_{2q-k-r}} (x,y,z) 
  \\
  & \quad \leq\sum_{y,z\in\N^{q-k}}\(\sum_{x\in\N^{k-r}}\sum_{\substack{(w_1,\ldots ,w_r,x_1,\ldots ,x_{k-r},y_1,\ldots ,y_{q-k})\in K\\(w_1,\ldots ,w_r,x_1,\ldots ,x_{k-r},z_1,\ldots ,z_{q-k})\in K}}\tilde f_q(w,x,y)\tilde f_q(w,x,z)\)^2{\bf 1}_{\tilde \Delta_{2q-2k}}  (y,z)
  \\
  & \quad \leq \left\| \(\tilde f_q\tilde\star_k^k\tilde f_q\){\bf 1}_{\tilde \Delta_{2q-2k}}
  \right\|_{l^2(\N)^{\circ (2q-2k)}},
\end{align*}
where the first inequality follows from the general inequality $\sum_{i\in I}a_i^2\leq\(\sum_{i\in I}a_i\)^2$, $a_i\geq0$, and the fact that $\text{sgn}\big(\tilde f_q(w,x,y)\tilde f_q(w,x,z)\big)$ is constant for fixed $y,z\in\N^{q-k}$. Thus, we  get
$$\left\| \(f_q\star_k^rf_q\){\bf 1}_{\Delta_{2q-k-r}}
\right\|_{L^2(\R_+)^{\circ (2q-k-r)}}\leq2^{2q}\(\E[X_1^4]\)^q\left\| \(\tilde f_q\tilde\star_k^k\tilde f_q\)
        {\bf 1}_{\tilde \Delta_{2q-2k}}
        \right\|_{l^2(\N)^{\circ (2q-2k)}}.$$
Analogously, for $r < k$ we obtain
$$\left\| f_q\star_k^rf_q\right\|_{L^2(\R_+)^{\circ (2q-k-r)}}\leq 2^{2q}\(\E[X_1^4]\)^q\left\| \(\tilde f_q\tilde\star_k^{k-1}\tilde f_q\)\right\|_{l^2(\N)^{\circ (2q-2k)}}.$$
Finally, applying this to (\ref{eq:aux3}), we may write
\begin{align*}
 &     d_W( I_q(f_q)
  \\
   &\leq C\(\E[X_1^4]\)^q\Bigg(
  \max_{k\in\{1,\ldots ,q-1\}}\left\| \(\tilde f_q\tilde\star_k^k\tilde f_q\)
      {\bf 1}_{\tilde \Delta_{2q-2k}}
      \right\|_{l^2(\N)^{\circ (2q-2k)}}
      +\max_{k\in\{1,\ldots ,q\}}\left\| \(\tilde f_q\tilde\star_k^{k-1}\tilde f_q\)\right\|_{l^2(\N)^{\circ (2q-2k)}}\Bigg),
\end{align*}
for some $C=C(q)$, and both maxima can be
calculated as in Theorem~6.2 of \cite{reichenbachs}.
\end{Proof}
Theorem~\ref{theoremk} extends the standard Berry-Esseen 
bound of Corollary~6.2 in  \cite{reichenbachs}
to general independent random sequences, in particular when
$K$ takes the form $K=\{1,\ldots ,n\}^q \cap \tilde\Delta_q$. 
Note also that the general result on random sequences
in Proposition~6.8 of \cite{nourdin3}
does not apply to the total variation or Wasserstein
distances. 
\section{Quadratic functionals} 
\label{s3.1}
\noindent
This section is devoted to double stochastic integrals, which are a special case of the multiple integrals discussed in Section~\ref{s3}.
We study them in a separate section because of many applications i.e. to
quadratic functionals. 
 Taking $n=2$ in Corollary~\ref{thm:dIZnabla} of Section~\ref{s3}, 
we get the following result.
 \begin{corollary}
  \label{thm:dIZnabla-2}
  Let $f_2\in L^2(\R_+^2)$ be a symmetric function
  satisfying \eqref{ass:int0}. 
  Assume that the functions 
  $$
  \widehat G^2_1f_2 ( y ) 
  = \frac{1}{2} \int_0^\infty | f_2(x , y ) |^2 dx 
\quad \mbox{and} \quad 
 \widehat G^2_2f_2 ( y , z ) 
 = \frac{1}{2} \int_0^\infty f_2(x, y ) f_2(x, z ) dx 
 $$
  belong to $L^2(\R_+)$ and $L^2(\R_+^2)$, respectively. Then we have
\begin{align*}
  & d_W( Z_n ,\mathcal N)
  \\
  & \leq \sqrt{ \(1-\frac14
  \Vert f_2\Vert_{L^2(\real_+^2)}^2 
  \)^2
  +\int_0^\infty \(\int_0^\infty | f_2(x,y) |^2\,dx\)^2\frac{dy}{2}
  +\int_0^\infty \int_0^\infty \(\int_0^\infty f_2(x,y)f_2(x,z)\,dx\)^2\frac{dy}{2} dz
  }\\
  &
  + \Vert f_2\Vert_{L^2(\real_+^2)} 
 \sqrt{ 
  \frac{1}{2} \int_0^\infty
  \hskip-0.2cm
  \(
  \int_0^\infty \hskip-0.2cm | f_2(x,y) |^2\,dx\)^2
  dy 
  + \Vert f_2\Vert_{L^4(\real_+^2)}^4 
  +\int_0^\infty \hskip-0.2cm \int_0^\infty \hskip-0.2cm
  \int_0^\infty \hskip-0.2cm
  | f_2(x,y)f_2(x,z) |^2dxdydz}.
\end{align*}
 \end{corollary}
 For example, when $f_2 \in C^1_{\Box}(\R_+^2) \cap L^2(\R_+^2 )$ is given by 
 \begin{equation}
   \label{f2} 
 f_2 (s,t) :=
\sum_{1 \leq k, l \leq n }
a_{k,l} f_1(s) f_1(t) {\bf 1}_{ ( 2k-2,2k]\times ( 2l-2,2l] } ( s,t), 
\qquad s , t \in \real_+,
\end{equation}
where $A=(a_{k,l})_{1\leq k,l\leq n}$ is a symmetric matrix with
vanishing diagonal
and such that $\sum_{1 \leq k , l \leq n} a_{k,l}^2 =1$, 
Corollary~\ref{thm:dIZnabla-2} yields the following result,
when $f_1$ is given by \eqref{djkldd}. 
\begin{corollary}
   \label{c01} 
  Given $(X_k)_{k\geq 1}$ a sequence of independent
  identically distributed random
 variables such that $E[X_k]=0$ and $E[X_k^2]=1$, $k\geq 1$,
 let $Q_n$ denote the normalized quadratic form 
 \begin{equation}
   \label{qn}
 Q_n : = \sum_{1 \leq k , l \leq n} a_{k,l} X_k X_l,
\end{equation}
 with $E[Q_n]=0$ and $E[Q_n^2]=1$, $n\geq 2$. We have 
 \begin{align}
   \nonumber
   & d_W( Q_n , \mathcal N)
    \\
    \nonumber
    & \leq 2 
   \sqrt{
    E[X_1^4] 
    \sum_{l=1}^n \( \sum_{k=1}^n a^2_{k,l} \)^2 
      +
    2 
             \sum_{1 \leq l , p \leq n }
\(\sum_{k=1}^n
      a_{k,l} a_{k,p}
      \)^2
      }
 +
   4 \sqrt{ 
\left( 3 
   E[X_1^4] 
     + ( E[X_1^4])^2 
\right) \sum_{k=1}^n 
\left( \sum_{l=1}^n a^2_{k,l} \right)^2 
}
.
\end{align} 
\end{corollary}
\begin{Proof}
 Writing $Q_n$ as 
$$
  Q_n : = 
 I_2 \big( f_2{\bf 1}_{[0,2n]\times [0,2n]} \big)
=
 \sum_{1 \leq k,l \leq n }
 a_{k,l}
 I_1 \big( f_1 {\bf 1}_{ ( 2k-2,2k] } \big) 
 I_1 \big( f_1 {\bf 1}_{ ( 2l-2,2l] } \big) 
,
$$ 
 $n\geq 2$, we have 
   \small 
   \begin{align*} 
     & d_W( Q_n , \mathcal N)
     \leq
 \frac{1}{\sqrt{2}}
  \sqrt{
    \int_0^{2n} \(\int_0^{2n} | f_2(x,y) |^2\,dx\)^2dy
    + \int_0^{2n} \int_0^{2n} \(\int_0^{2n} f_2(x,y)f_2(x,z)\,dx\)^2dydz }
  \\
  &+
  2 
  \sqrt{ 
  \frac{1}{2} \int_0^{2n}
  \hskip-0.1cm
  \left|
  \int_0^{2n} \hskip-0.2cm | f_2(x,y) |^2\,dx\right|^2
  dy 
  +\int_0^{2n} \hskip-0.2cm \int_0^{2n}
  \hskip-0.2cm
  | f_2(x,y) |^4dxdy
  +\int_0^{2n} \hskip-0.2cm \int_0^{2n} \hskip-0.2cm
  \int_0^{2n} \hskip-0.2cm
  | f_2(x,y)f_2(x,z) |^2dxdydz}
\\
  & = 
   \frac{1}{\sqrt{2}} 
   \int_0^2 |f_1(x)|^2 \,dx
   \sqrt{
    \int_0^2 |f_1(y)|^4 dy
    \sum_{l=1}^n \( \sum_{k=1}^n a^2_{k,l} \)^2 
      +
      \left( \int_0^2 |f_1(x)|^2 \,dx \right)^2 
            \sum_{1 \leq l , p \leq n }
\(\sum_{k=1}^n
      a_{k,l} a_{k,p}
      \)^2 }
  \\
  &+
   2 \sqrt{ 
 \left|
  \int_0^2 \hskip-0.2cm |f_1(x)|^2\,dx \right|^2 
     \int_0^2 \hskip-0.1cm |f_1(y)|^4 dy 
     \left(
 \frac{1}{2}     \sum_{l=1}^n 
      \left(
   \sum_{k=1}^n a_{k,l}^2
   \right)^2
   +
  \sum_{1 \leq l , p \leq n }
\sum_{k=1}^n
 a^2_{k,l} a^2_{k,p}
\right) 
 +\left| \int_0^2 \hskip-0.2cm |f_1(x)|^4 dx \right|^2 
\sum_{1 \leq k, l \leq n } a^4_{k,l} 
}
  \\
  & = 
 2 E[X_1^2] 
   \sqrt{
    E[X_1^4] 
    \sum_{l=1}^n \( \sum_{k=1}^n a^2_{k,l} \)^2 
      +
    2 ( E[X_1^2] )^2 
             \sum_{1 \leq l , p \leq n }
\(\sum_{k=1}^n
      a_{k,l} a_{k,p}
      \)^2
      }
  \\
  &+
   4 \sqrt{ 
 3 ( E[X_1^2])^2 
   E[X_1^4] 
      \sum_{l=1}^n 
      \left(
   \sum_{k=1}^n a_{k,l}^2
   \right)^2
 + ( E[X_1^4])^2 
\sum_{k=1}^n 
\left( \sum_{l=1}^n a^2_{k,l} \right)^2 
}
, 
\end{align*}

   \normalsize

   \noindent
   where we used the relation 
$$
\int_0^2 | f_1(x) |^4dx = 2 \int_0^1 | F^{-1} (y) |^4 dy = 2 E[X_1^4].
$$
\end{Proof} 
Bounds of that type have been already studied in the literature,
see e.g. \cite{rotar} and \cite{goetze}.
They are usually presented by means of the expression
$$
L_n^2 : = \max_{1\leq k \leq n} \sum_{l=1}^n a_{k,l}^2.
$$
Following this convention we can apply  the bound
of Corollary~\ref{c01} to obtain 
\begin{equation}
\label{jdk2} 
  d_W ( Q_n , \mathcal N)\leq
  2 \sqrt{n} L_n^2 \left(
  \sqrt{
     E[X_1^4] 
      +
    \frac{2}{nL_n^4} \sum_{1 \leq l , p \leq n }
\(\sum_{k=1}^n
      a_{k,l} a_{k,p}
      \)^2
      } +
   2 \sqrt{ 
 3 E[X_1^4] 
 + ( E[X_1^4])^2 
}\right). 
\end{equation} 
 Note that the constants in the above bound are explicit. 
For example, when
   $$
   Q_n = \frac{2}{\sqrt{n}} \sum_{k=1}^n X_{2k-1}X_{2k},
   $$
   we have $L_n^2 = 1/n$ and
   $$
   \sum_{1 \leq l , p \leq 2n }
\(\sum_{k=1}^{2n} 
      a_{k,l} a_{k,p}
      \)^2 = \frac{1}{n}, 
      $$
 hence \eqref{jdk2} recovers the known convergence rate 
$$
 d_W( Q_n , \mathcal N)
 \leq
  \sqrt{\frac{2}{n}} \left(
  \sqrt{
    1 + E[X_1^4] 
      } +
   2 \sqrt{ 
 3 E[X_1^4] 
 + ( E[X_1^4])^2 
   }\right)
   \leq
  \frac{8E[X_1^4]}{\sqrt{n}}
   , 
$$    
  cf. pages~1074-1075 of \cite{goetze},
  with an explicit constant
  depending on $E[X_1^4]$ instead of $\sqrt {E[|X_1|^3]}$.
 On the other hand, Corollary~\ref{thm:dIZD} applied with $n=2$
 gives the following result. 
\begin{corollary}\label{thm:dIZD-2}
  For any $f_2 \in \mathcal C^1_{\Box}(\R_+^2) \cap L^2(\R_+^2)$
  satisfying \eqref{ass:int0}, we have 
  \begin{align}
    \nonumber 
&    d_W( I_2(f_2),\mathcal N)
  \\
    \nonumber 
  & \leq \left|1-\frac14
  \Vert f_2 \Vert_{L^2(\real_+^2)}^2\right|
+\Bigg\{2\int_0^\infty \(\int_0^x \partial_1f_2(x,y)\int_0^{x}f_2(t,y)dt\,dy\)^2dx\\
    \nonumber 
& +4\int_0^\infty\int_0^x\(\partial_1f_2(x,y)\int_0^{x}f_2(t,y)dt\)^2dydx
 +4\int_0^\infty \int_0^\infty \int_0^x \(\partial_1f_2(x,y)\int_0^{x}f_2(t,z)dt\)^2dzdxdy\\
    \nonumber 
 & +2\int_0^\infty\(\int_y^\infty | f_2(x,y) |^2 dx \)^2dy+2
 \int_0^\infty \int_0^\infty \(\int_z^\infty f_2(x,y)f_2(x,z)dx \)^2dydz\\
 \label{nnb}
& -4\sum_{i=0}^\infty\(\int_{2i}^{2i+2}\int_0^x | f_2(x,y) |^2dydx\)^2\Bigg\}^{1/2}.
\end{align}
The bound for $d_{TV}( Q_n , \mathcal N)$ is twice as large
as \eqref{nnb}. 
    \end{corollary}
\begin{Proof}
 We apply Corollary~\ref{thm:dIZD} with 
$$
H_0 ( x ) = 
 \frac12 \int_0^x \partial_1f_2 (x,y) \int_0^x f_2 (t,y) dt dy,
 \quad
H_1 ( x, y ) = 
\mathbf1_{\{y<x \}} \partial_1f_2 (x,y) \int_0^x f_2 (t,y) dt 
,
$$
$$
H_2 ( x, y , z ) = 
\mathbf1_{\{z<x \}}
\frac{1}{2} \partial_1f_2 (x,y) \int_0^x f_2 (t,z) dt
+
\mathbf1_{\{y<x \}} \frac{1}{2} \partial_1f_2 (x,z) \int_0^x f_2 (t,y) dt
,
 \quad
 x,y,z\in \real_+,
$$
 and
 $$
 J_1(s,y) = | f_2(s,y) |^2 \mathbf1_{\{y<s\}}, 
  \quad 
  J_2(s,y,z) = \frac{1}{2} f_2(s,y) f_2(s,z) \mathbf1_{\{z<s\}}
  +
  \frac{1}{2} f_2(s,y) f_2(s,z) \mathbf1_{\{y<s\}}, 
$$
  $$
  \hat{J}_1(z_1) = \frac{1}{2} \int_{z_1}^\infty ( f_2(s,z))^2 ds,
  \quad 
  \hat{J}_2(z_1,z_2)
  = \frac{1}{4} \int_{z_2}^\infty f_2(s,z_1) f_2(s,z_2) ds 
  + \frac{1}{4} \int_{z_1}^\infty f_2(s,z_1) f_2(s,z_2) ds. 
$$
\end{Proof} 
 When $f_2 \in C^1_{\Box}(\R_+^2) \cap L^2(\R_+^2 )$ is given by
\eqref{f2}, Corollary~\ref{thm:dIZD-2} shows the following
bound on quadratic functionals.
\begin{corollary}
   \label{c01.1} 
  Given $(X_k)_{k\geq 1}$ a sequence of independent
  identically distributed random
 variables such that $E[X_k]=0$ and $E[X_k^2]=1$, $k\geq 1$,
 the normalized quadratic form $Q_n$ defined in
 \eqref{qn} satisfies 
\begin{equation} 
  \label{qn2}
  d_W( Q_n , \mathcal N)
   \leq 
 4 \sqrt{
  E [ ( \varphi_{X_k} (X_k) )^2 ]
 ( 2 + E[X_1^4] )  L_n^2
     + 
     2
     \hskip-0.2cm
     \sum_{1 \leq q , l \leq n }
     \hskip-0.1cm
     \(
    \sum_{k=1}^n 
               a_{k,q}a_{k,l}  
               \)^2
                    \hskip-0.2cm
                    -  \sum_{k=1}^n \(  \sum_{l=1}^{k-1}  a^2_{k,l}\)^2
}.
\end{equation} 
The bound for $d_{TV}( Q_n , \mathcal N)$ is twice as large
as \eqref{qn2}. 
 \end{corollary}
 \begin{Proof} 
 By Corollary~\ref{thm:dIZD-2}, we have 
 \begin{align*}
   &    d_W( I_2(f_2),\mathcal N) \leq
   \sqrt{
     2 I_1 + 4 I_2 + 4 I_3 + 2 I_4 + 2 I_5 - 4 I_6
  }, 
 \end{align*}
 where
 \begin{align*}
 I_1&=\int_0^\infty \(\int_0^x
\sum_{1 \leq k, l , p, q \leq n }
a_{k,l} f'_1(x) f_1(y) {\bf 1}_{ ( 2k-2,2k]\times ( 2l-2,2l] } ( x,y)
    \right.\\
    &\hspace{30mm}\left.\times\int_0^{x}
    a_{p,q} f_1(t) f_1(y) {\bf 1}_{ ( 2p-2,2p]\times ( 2q-2,2q] } ( t,y)
dt\,dy\)^2dx\\
&=\int_0^\infty 
  \sum_{k=1}^n {\bf 1}_{ ( 2k-2,2k] } ( x ) 
\(
\sum_{1 \leq p , l \leq k }
    a_{k,l} a_{p,l} 
      \int_0^x |f_1(y)|^2 {\bf 1}_{ ( 2l-2,2l] } (y)
          dy
    f'_1(x)\right.\\
    &\hspace{30mm}\left.\times\int_0^{x} f_1(t) {\bf 1}_{ ( 2p-2,2p] } ( t)
          dt\,
          \)^2dx\\
          &=\int_0^\infty 
  \sum_{k=1}^n {\bf 1}_{ ( 2k-2,2k] } ( x ) 
\(
\sum_{1 \leq  l \leq k }
    (a_{k,l})^2  
      \int_0^x |f_1(y)|^2 {\bf 1}_{ ( 2l-2,2l] } (y)
          dy
    f'_1(x)\right.\\
    &\hspace{30mm}\left.\times\int_0^{x} f_1(t) {\bf 1}_{ ( 2k-2,2k] } ( t)
          dt\,
          \)^2dx\\
          &=\( \int_0^2 |f_1(y)|^2 dy \)^2 \int_0^2 
\(   f'_1(x) \int_0^x f_1(t) dt \)^2dx  
  \sum_{k=1}^n 
\(
\sum_{l=1 }^{k-1}
    (a_{k,l})^2
          \)^2\\
          &\leq 8 (E[X_1^2])^2
  E [ ( \varphi_{X_k} (X_k) )^2 ]
    \sum_{k=1}^n 
\(
\sum_{l=1 }^{k-1} (a_{k,l})^2  \)^2,
 \end{align*}
   \begin{align*}
 I_2&=\int_0^\infty\int_0^x\(
\sum_{1 \leq k, l \leq n }
a_{k,l} f'_1(x) f_1(y) {\bf 1}_{ ( 2k-2,2k]\times ( 2l-2,2l] } ( x,y)
    \right.\\
    &\hspace{30mm}\left.\times\int_0^{x}
    \sum_{1 \leq p, q \leq n }
    a_{p,q} f_1(t) f_1(y) {\bf 1}_{ ( 2p-2,2p]\times ( 2q-2,2q] } ( t,y)
        dt\)^2dydx\\
        &=\int_0^\infty\int_0^x\(
\sum_{1 \leq k, l \leq n }
a_{k,l} f'_1(x) f_1(y) {\bf 1}_{ ( 2k-2,2k]\times ( 2l-2,2l] } ( x,y)
    \right.\\
    &\hspace{30mm}\left.\times\int_0^{x}
    a_{k,l} f_1(t) f_1(y) {\bf 1}_{ ( 2k-2,2k] } ( t)
        dt\)^2dydx\\
        &=\sum_{1 \leq k, l \leq n }
          a^4_{k,l}
          \int_0^\infty 
          {\bf 1}_{ ( 2k-2,2k] } ( x)
          \(
    f'_1(x) \int_0^{x}
 f_1(t)
    {\bf 1}_{ ( 2k-2,2k] } ( t)
        dt\)^2        \\
        &\hspace{30mm}\times        \int_0^x |f_1(y)|^4 
{\bf 1}_{ ( 2l-2,2l] } ( y)
dydx\\
&\leq \sum_{1 \leq k, l \leq n }
          a^4_{k,l}
          \int_0^2 
          \( f'_1(x) \int_0^x 
 f_1(t)        dt \)^2 dx     \int_0^2 |f_1(y)|^4 
dy\\
&\leq 4 E[X_1^4] 
E [ ( \varphi_{X_k} (X_k) )^2 ]
\sum_{1 \leq k, l \leq n }
          a^4_{k,l},
 \end{align*}
  \begin{align*}
 I_3&=\int_0^\infty \int_0^\infty \int_0^x
        \(
        \sum_{1 \leq k, l  \leq n }
a_{k,l} f'_1(x) f_1(y) {\bf 1}_{ ( 2k-2,2k]\times ( 2l-2,2l] } ( x,y)
    \right.\\
    &\hspace{30mm}\left.\times\int_0^{x}
        \sum_{1 \leq p, q \leq n }
    a_{p,q} f_1(t) f_1(z) {\bf 1}_{ ( 2p-2,2p]\times ( 2q-2,2q] } ( t,z)
        dt\)^2dzdxdy\\
        &=\int_0^\infty \int_0^\infty \int_0^x
        \(
        \sum_{1 \leq k, l  \leq n }
a_{k,l} f'_1(x) f_1(y) {\bf 1}_{ ( 2k-2,2k]\times ( 2l-2,2l] } ( x,y)
    \right.\\
    &\hspace{30mm}\left.\times\int_0^{x}
        \sum_{1 \leq q \leq n }
    a_{k,q} f_1(t) f_1(z) {\bf 1}_{ ( 2k-2,2k]\times ( 2q-2,2q] } ( t,z)
        dt\)^2dzdxdy\\
        &=\sum_{1 \leq k, l , q \leq n }
        \int_0^\infty \int_0^\infty \int_0^x
        {\bf 1}_{ ( 2k-2,2k]\times ( 2l-2,2l] } ( x,y)
    {\bf 1}_{ ( 2q-2,2q] } ( z)
  |f_1(y)|^2 |f_1(z)|^2 
        a^2_{k,l} \\
       &\hspace{30mm}\times \(
    a_{k,q} f'_1(x) \int_0^{x}
   f_1(t) 
    {\bf 1}_{ ( 2p-2,2p] } ( t)
        dt\)^2dzdxdy\\
        &\leq  \sum_{1 \leq k, l , q \leq n }
         \( \int_0^2 |f_1(y)|^2 dy \)^2  
        \( \int_0^2 f'_1(x) \int_0^x f_1(y) dy \)^2 
        a^2_{k,l}   a^2_{k,q} \\
    &\leq  4 
        (E[X_1^2])^2 
        E [ ( \varphi_{X_k} (X_k) )^2 ]
        \sum_{1 \leq k, l , q \leq n }
        a^2_{k,l}a^2_{k,q} ,
 \end{align*}
 hence
  $$
  2 I_1 + 4 I_2 + 4 I_3
   \leq 
  16
  E [ ( \varphi_{X_k} (X_k) )^2 ]
  \left(
    \sum_{k=1}^n 
\(
\sum_{l=1 }^{k-1} (a_{k,l})^2  \)^2
+ 
 E[X_1^4] \sum_{1 \leq k, l \leq n } a^4_{k,l}
 + 
  \sum_{1 \leq k, l , q \leq n }
        a^2_{k,l}a^2_{k,q}  \right).
 $$ 
    On the other hand, we have 
\begin{align*}
 I_4&=\frac12\int_0^\infty\(\int_y^\infty \left|
                \sum_{1 \leq k, l \leq n }
a_{k,l} f_1(x) f_1(y) {\bf 1}_{ ( 2k-2,2k]\times ( 2l-2,2l] } ( x,y)
    \right|^2 dx \)^2dy\\
    &= \frac12
        \sum_{1 \leq k, l \leq n } a^4_{k,l} 
        \int_0^\infty                 {\bf 1}_{ ( 2l-2,2l] } ( y)
          |f_1(y)|^4 \( 
          \int_y^\infty |f_1(x)|^2 
                {\bf 1}_{ ( 2k-2,2k] } ( x)
    dx \)^2 dy\\
    &\leq \frac12
        \sum_{1 \leq k, l \leq n } a^4_{k,l} 
        \( \int_0^2 |f_1(y)|^2 dy \)^4
        \leq 8 \( E[X_1^2] \)^4  \sum_{1 \leq k, l \leq n } a^4_{k,l} ,
 \end{align*}
  \begin{align*}
 I_5&=\int_0^\infty \int_0^\infty \(\int_z^\infty
                    \sum_{1 \leq k, l , p, q \leq n }
a_{k,l} f_1(x) f_1(y) {\bf 1}_{ ( 2k-2,2k]\times ( 2l-2,2l] } ( x,y)\right.\\
&\hspace{40mm}\left.\times a_{p,q} f_1(x) f_1(z) {\bf 1}_{ ( 2k-2,2k]\times ( 2l-2,2l] } ( x,z)
    dx \)^2dydz\\
    &= \int_0^\infty \int_0^\infty
    |f_1(y)|^2 |f_1(z)|^2 
    \sum_{1 \leq q , l \leq n }
                       {\bf 1}_{ ( 2l-2,2l] } ( y)
                  {\bf 1}_{ ( 2q-2,2q] } ( z)\\
                  &\hspace{40mm} 
                  \times\( \sum_{k=1}^n 
               a_{k,q}a_{k,l}                    \int_z^\infty
                    |f_1(x)|^2 
                    {\bf 1}_{ ( 2k-2,2k] } ( x)
    dx \)^2dydz\\
    &\leq \( \int_0^2
    |f_1(y)|^2 dy \)^4 
    \sum_{1 \leq q , l \leq n }
    \(
    \sum_{k=1}^n 
               a_{k,q}a_{k,l}  
               \)^2 \leq 16 
    \( E[X_1^2] \)^4 
    \sum_{1 \leq q , l \leq n }
    \(
    \sum_{k=1}^n 
               a_{k,q}a_{k,l}  
               \)^2. 
 \end{align*}
 and
\begin{align*}
I_6&=\sum_{i=1}^n \(\int_{2i-2}^{2i}\int_0^x
\left|
                \sum_{1 \leq k, l \leq n }
a_{k,l} f_1(x) f_1(y) {\bf 1}_{ ( 2k-2,2k]\times ( 2l-2,2l] } ( x,y)
    \right|^2dydx\)^2\\
    &=\sum_{k=1}^n \(
     \sum_{l=1}^n 
    a^2_{k,l}
    \int_{2k-2}^{2k}
      |f_1(x)|^2 \int_0^x
        |f_1(y)|^2 
            {\bf 1}_{ ( 2l-2,2l] } ( y)
              dydx
              \)^2\\
              &= \sum_{k=1}^n \( 
                   \sum_{l=1}^{k-1} 
    a^2_{k,l}
    \int_{2k-2}^{2k}
      |f_1(x)|^2 \int_{2l-2}^{2l}
        |f_1(y)|^2 
              dydx
\)^2\\
&= 16 
    \( E[X_1^2] \)^4 \sum_{k=1}^n \( 
                   \sum_{l=1}^{k-1} 
    a^2_{k,l}
\)^2.
\end{align*}
Hence we have 
\begin{align*}
  2I_4+2I_5-4I_6&\leq 32
  \(\sum_{1 \leq k, l \leq n } a^4_{k,l}+\sum_{1 \leq q , l \leq n }
    \(\sum_{k=1}^n            a_{k,q}a_{k,l}  \)^2- 2 \sum_{k=1}^n \(  \sum_{l=1}^{k-1}  a^2_{k,l}\)^2\)\\
    &\leq 32 \(\sum_{1 \leq q , l \leq n }
    \(\sum_{k=1}^n            a_{k,q}a_{k,l}  \)^2-  \sum_{k=1}^n \(  \sum_{l=1}^{k-1}  a^2_{k,l}\)^2\), 
\end{align*}
 and combining the above bounds gives us \eqref{qn2}. 
\end{Proof} 
 When $(X_k)_{k\geq 1}$ is a sequence of independent
 gamma identically distributed normalized random
 variables we have $E [ ( \varphi_{X_k} (X_k) )^2 ] =2$, 
 and \eqref{qn2} yields 
$$ 
 d_W( Q_n , \mathcal N)
 \leq 
 4 \sqrt{
2 ( 2 + E[X_1^4] )L_n^2 
     + 
     2
     \hskip-0.2cm
     \sum_{1 \leq q , l \leq n }
     \hskip-0.1cm
     \(
    \sum_{k=1}^n 
               a_{k,q}a_{k,l}  
               \)^2
                    \hskip-0.2cm
                    -  \sum_{k=1}^n \(  \sum_{l=1}^{k-1}  a^2_{k,l}\)^2
}.
$$
  A similar expression can be obtained from \eqref{qn3} in the
  beta case.
\section{Appendix - multiplication formula} 
\label{s5}
We now formulate and prove the multiplication formula which
is used in the proof of Corollary~\ref{thm:dIZnabla}. 
\begin{theorem}
  \label{thmm}
  (Multiplication formula).
  Let  $f_n$, $g_m$ satisfy \eqref{ass:int0}  and $f_n\star_{k}^{i}g_m\in L^2(\R_+^{m+n-k-i})$ for all $0\leq i \leq k\leq m\wedge n$. Then we have
\begin{equation}
\nonumber 
 I_n(f_n) I_m(g_m)=\sum_{k=0}^{m\wedge n}k!{{m}\choose{k}}{{n}\choose{k}}\sum_{i=0}^k{{k}\choose{i}} I_{m+n-k-i}\(f_n\hskip0.1cm \widetilde{\star}_{k}^{i}g_m\).
\end{equation}
\end{theorem}
\begin{Proof} Without loss of generality we consider only $n\geq m$. We use mathematical induction with respect to $m$, if $m<n$, and with respect to $n$, if $m=n$.  The formula is clearly valid for $n\geq 0$ and $m=0$.
    Let us assume that the formula is valid for the following pairs of indices: $(n,m-1)$, $(n-1,m)$ and $(n-1,m-1)$.
By \eqref{eq:Ipsi} we get
\begin{align*}
 \big( I_n(f_n&) I_m(g_m)\big)\circ \Psi_t =S_1(t)+S_2(t)+S_3(t),
 \end{align*}
 where 
 $$
 \left\{
 \begin{array}{l} 
    S_1(t) =mn I_{n-1}(f_n(t,*)) I_{m-1}(g_m(t,*))+m I_{n}(f_n) I_{m-1}(g_m(t,*))+n I_{n-1}(f_n(t,*)) I_{m}(g_m),
   \\ \\
 S_2(t) =-\big(n I_{n-1}(f_n(v,*)) I_{m}(g_m) \circ \Psi_t +m I_{n}(f_n)  \circ \Psi_t I_{m-1}(g_m(v,*))\big)_{\mid v=2\lfloor t/2\rfloor+1+U_{\lfloor t/2\rfloor}},\\ \\
 S_3(t) = I_n(f_n) I_m(g_m)-mn I_{n-1}(f_n(v,*)) I_{m-1}(g_m(v,*))_{\mid v=2\lfloor t/2\rfloor+1+U_{\lfloor t/2\rfloor}}.
 \end{array}
 \right.
 $$ 
We note that by \eqref{ass:int0} we have $\E\big[S_2(t) \mid \tilde{{\cal F}}_t\big]=0$. Additionally, the function $s\longmapsto \E\big[S_3(t) \mid \tilde{{\cal F}}_s \big]$ is constant for $s \in[2\lfloor t/2\rfloor,2\lfloor t/2\rfloor+2)$
  which, combined with \eqref{eq:C-O3}, implies
\begin{eqnarray*}
\int_0^\infty \E\big[\nabla_t( I_n(f_n) I_m(f_m)) \mid \tilde{{\cal F}}_t \big]d(Y_t-t/2)
&=&\int_0^\infty\E\big[\big( I_n(f_n) I_m(g_m)\big) \circ \Psi_t \mid \tilde{{\cal F}}_t\big]d(Y_t-t/2)\\
&=&\int_0^\infty\E\big[S_1(t) \mid \tilde{{\cal F}}_t\big]d(Y_t-t/2).
\end{eqnarray*}
Then, by the induction hypothesis and renumeration in the first sum below, we get
\begin{align*}
  &\int_0^\infty\E\big[
    \nabla_t( I_n(f_n) I_m(f_m)) \mid \tilde{{\cal F}}_t\big]
  d(Y_t-t/2)\\
&=\int_0^\infty E \left[ nm\sum_{k=0}^{(m\wedge n)-1}k!{{m-1}\choose{k}}{{n-1}\choose{k}}\sum_{i=0}^k{{k}\choose{i}} I_{m+n-2-k-i}\(f_n(t,*)\star_{k}^{i}g_m(t,*)\)\right.\\
&\qquad \left.+m\sum_{k=0}^{(m-1)\wedge n}k!{{m-1}\choose{k}}{{n}\choose{k}}\sum_{i=0}^k{{k}\choose{i}} I_{m+n-1-k-i}\(f_n\star_{k}^{i}g_m(t,*)\)\right.\\
  &\qquad \left.+n\sum_{k=0}^{m\wedge (n-1)}k!{{m}\choose{k}}{{n-1}\choose{k}}\sum_{i=0}^k{{k}\choose{i}} I_{m+n-1-k-i}\(f_n(t,*)\star_{k}^{i}g_m\)\Big| \tilde{{\cal F}}_t \right]
d(Y_t-t/2)\\
&=\int_0^\infty
E \left[ \sum_{k=1}^{(m\wedge n)}k!{{m}\choose{k}}{{n}\choose{k}}\sum_{i=0}^{k-1}(k-i){{k}\choose{i}} I_{m+n-1-k-i}\(f_n(t,*)\star_{k-1}^{i}g_m(t,*)\)\right.\\
&\qquad \left.+\sum_{k=0}^{(m-1)\wedge n}k!(m-k){{m}\choose{k}}{{n}\choose{k}}\sum_{i=0}^k{{k}\choose{i}} I_{m+n-1-k-i}\(f_n\star_{k}^{i}g_m(t,*)\)\right.\\
&\qquad \left.+\sum_{k=0}^{m\wedge (n-1)}k!(n-k){{m}\choose{k}}{{n}\choose{k}}\sum_{i=0}^k{{k}\choose{i}} I_{m+n-1-k-i}\(f_n(t,*)\star_{k}^{i}g_m\) \Big| \tilde{{\cal F}}_t \right] d(Y_t-t/2).
\end{align*}
Thus, the formula \eqref{eq:COI} gives us for $m\neq n$
\begin{align*}
  \int_0^\infty&\E\big[
    \nabla_t( I_n(f_n) I_m(f_m)) \mid \tilde{{\cal F}}_t \big]
  d(Y_t-t/2)\\
&= \sum_{k=1}^{(m\wedge n)}k!{{m}\choose{k}}{{n}\choose{k}}\sum_{i=0}^{k-1}\frac{k-i}{m+n-k-i}{{k}\choose{i}} I_{m+n-k-i}\(f_n\star_{k}^{i}g_m\)\\
&\quad \left.+\sum_{k=0}^{(m-1)\wedge n}k!{{m}\choose{k}}{{n}\choose{k}}\sum_{i=0}^k\frac{m-k}{m+n-k-i}{{k}\choose{i}} I_{m+n-k-i}\(f_n\star_{k}^{i}g_m\)\right.\\
&\quad +\sum_{k=0}^{m\wedge (n-1)}k!{{m}\choose{k}}{{n}\choose{k}}\sum_{i=0}^k\frac{n-k}{m+n-k-i}{{k}\choose{i}} I_{m+n-k-i}\(f_n\star_{k}^{i}g_m\). 
\end{align*}
If $m=n$, then the term $I_0\(f_n\star_{n}^{n}g_n\)$ does not appear in the above sums. Nevertheless, since $I_0\(f_n\star_{n}^{n}g_n\)=\langle f_n,g_n\rangle_{L^2(\R^n_+,dx/2)}$, the assertion of the theorem follows from \eqref{eq:EI^2} and Clark-Ocone formula \eqref{eq:C-O2}.
\end{Proof}

\footnotesize

\def\cprime{$'$} \def\polhk#1{\setbox0=\hbox{#1}{\ooalign{\hidewidth
  \lower1.5ex\hbox{`}\hidewidth\crcr\unhbox0}}}
  \def\polhk#1{\setbox0=\hbox{#1}{\ooalign{\hidewidth
  \lower1.5ex\hbox{`}\hidewidth\crcr\unhbox0}}} \def\cprime{$'$}


\begin{thebibliography}{10}

\bibitem{blei}
R.~Blei and S.~Janson.
\newblock Rademacher chaos: tail estimates versus limit theorems.
\newblock {\em Ark. Mat.}, 42(1):13--29, 2004.

\bibitem{goldstein}
L.~Goldstein.
\newblock Bounds on the constant in the mean central limit theorem.
\newblock {\em Ann. Probab.}, 38(4):1672--1689, 2010.

\bibitem{goetze}
F.~G\"otze and A.N. Tikhomirov.
\newblock Asymptotic distribution of quadratic forms.
\newblock {\em Ann. Probab.}, 27(2):1072--1098, 1999.

\bibitem{reichenbachs}
K.~Krokowski, A.~Reichenbachs, and C.~Thaele.
\newblock Berry-{E}sseen bounds and multivariate limit theorems for functionals
  of {R}ademacher sequences.
\newblock {\em Ann. Inst. Henri Poincar\'e Probab. Stat.}, 52(2):763--803,
  2016.

\bibitem{krokowski}
K.~Krokowski, A.~Reichenbachs, and C.~Thaele.
\newblock Discrete {M}alliavin-{S}tein method: {B}erry-{E}sseen bounds for
  random graphs and percolation.
\newblock {\em Ann. Probab.}, 45(2):1071--1109, 2017.

\bibitem{lnp}
M.~Ledoux, I.~Nourdin, and G.~Peccati.
\newblock Stein's method, logarithmic {S}obolev and transport inequalities.
\newblock {\em Geom. Funct. Anal.}, 25:256--306, 2015.

\bibitem{ley}
C.~Ley, G.~Reinert, and Y.~Swan.
\newblock Stein's method for comparison of univariate distributions.
\newblock {\em Probab. Surv.}, 14:1--52 (electronic), 2017.

\bibitem{nourdin}
I.~Nourdin.
\newblock Lectures on {G}aussian approximations with {M}alliavin calculus.
\newblock In {\em S\'eminaire de {P}robabilit\'es {XLV}}, volume 2078 of {\em
  Lecture Notes in Math.}, pages 3--89. Springer, 2013.

\bibitem{nourdinpeccati}
I.~Nourdin and G.~Peccati.
\newblock Stein's method on {W}iener chaos.
\newblock {\em Probab. Theory Related Fields}, 145(1-2):75--118, 2009.

\bibitem{nourdin3}
I.~Nourdin, G.~Peccati, and G.~Reinert.
\newblock Stein's method and stochastic analysis of {R}ademacher functionals.
\newblock {\em Electron. J. Probab.}, 15:no. 55, 1703--1742, 2010.

\bibitem{viens}
I.~Nourdin and F.G. Viens.
\newblock Density formula and concentration inequalities with {M}alliavin
  calculus.
\newblock {\em Electron. J. Probab.}, 14:2287--2309, 2009.

\bibitem{utzet2}
G.~Peccati, J.~L. Sol{\'e}, M.~S. Taqqu, and F.~Utzet.
\newblock Stein's method and normal approximation of {P}oisson functionals.
\newblock {\em Ann. Probab.}, 38(2):443--478, 2010.

\bibitem{thale}
G.~Peccati and C.~Th{\"a}le.
\newblock Gamma limits and {U}-statistics on the {P}oisson space.
\newblock {\em ALEA Lat. Am. J. Probab. Math. Stat.}, 10(1):525--560, 2013.

\bibitem{prebub}
N.~Privault.
\newblock Calcul des variations stochastique pour la mesure de densit{\'e}
  uniforme.
\newblock {\em {P}otential {A}nalysis}, 7(2):577--601, 1997.

\bibitem{girunif}
N.~Privault.
\newblock Absolute continuity in infinite dimensions and anticipating
  stochastic calculus.
\newblock {\em {P}otential {A}nalysis}, 8(4):325--343, 1998.

\bibitem{qtmunif}
N.~Privault.
\newblock Quantum stochastic calculus for the uniform measure and {B}oolean
  convolution.
\newblock In {\em S\'eminaire de Probabilit\'es, XXXV}, pages 28--47. Springer,
  Berlin, 2001.

\bibitem{privaultbk2}
N.~Privault.
\newblock {\em Stochastic analysis in discrete and continuous settings with
  normal martingales}, volume 1982 of {\em Lecture Notes in Mathematics}.
\newblock Springer-Verlag, Berlin, 2009.

\bibitem{privaulttorrisi3}
N.~Privault and G.L. Torrisi.
\newblock Probability approximation by {C}lark-{O}cone covariance
  representation.
\newblock {\em Electron. J. Probab.}, 18:1--25, 2013.

\bibitem{privaulttorrisi4}
N.~Privault and G.L. Torrisi.
\newblock The {S}tein and {C}hen-{S}tein methods for functionals of
  non-symmetric {B}ernoulli processes.
\newblock {\em ALEA Lat. Am. J. Probab. Math. Stat.}, 12:309--356, 2015.

\bibitem{rotar}
V.I. Rotar' and T.~L. Shervashidze.
\newblock Some estimates for distributions of quadratic forms. ({R}ussian).
\newblock {\em Teor. Verojatnost. i Primenen.}, 30(3):549--554, 1985.

\bibitem{stein}
C.~Stein.
\newblock Estimation of the mean of a multivariate normal distribution.
\newblock {\em Ann. Stat.}, 9(6):1135--1151, 1981.

\bibitem{surgailis}
D.~Surgailis.
\newblock On multiple {P}oisson stochastic integrals and associated {M}arkov
  semi-groups.
\newblock {\em Probability and Mathematical Statistics}, {3}:217--239, 1984.

\bibitem{surgailisclt}
D.~Surgailis.
\newblock Non-{CLT}s: {$U$}-statistics, multinomial formula and approximations
  of multiple {I}t\^o-{W}iener integrals.
\newblock In {\em Theory and applications of long-range dependence}, pages
  129--142. Birkh\"auser Boston, Boston, MA, 2003.

\end{thebibliography}
\end{document}